\documentclass[11pt]{article}
\usepackage{amsmath,amsfonts}
\usepackage{algorithmic}
\usepackage{algorithm}
\usepackage{array}
\usepackage[caption=false,font=normalsize,labelfont=sf,textfont=sf]{subfig}
\usepackage{textcomp}
\usepackage{stfloats}
\usepackage{url}
\usepackage{verbatim}
\usepackage{graphicx}
\usepackage{cite}
\usepackage{cuted} 

\usepackage{color}
\usepackage[normalem]{ulem}

\def\0{\bf \0}
\def\A{{\bf A}}
\def\B{{\bf B}}

\def\I{{\bf I}}
\def\J{{\bf J}}
\def\K{{\bf K}}
\def\L{{\bf L}}
\def\M{{\bf M}}

\def\0{{\bf 0}}

\def\Q{{\bf Q}}

\def\T{{\bf T}}

\def\V{{\bf V}}

\def\X{{\bf X}}

\def\a{{\bf a}}
\def\b{{\bf b}}

\def\f{{\bf f}}

\def\h{{\bf h}}

\def\p{{\bf p}}

\def\r{{\bf r}}

\def\u{{\bf u}}
\def\v{{\bf v}}

\def\x{{\bf x}}
\def\y{{\bf y}}
\def\z{{\bf z}}

\def\T{{\rm T}}
\def\diag{{\rm diag}}

\newtheorem{remark}{Remark}[section]
\begin{document}

\title{A Systematic Methodology for Modeling and Attitude Control of Multi-body Space Telescopes}

\author{Yaguang Yang, William Bentz, Lia Lewis
\thanks{Aerospace Engineers, Mission Engineering and Systems Analysis Division, Goddard Space Flight Center, NASA, 8800 Greenbelt Rd, Greenbelt, 20771 MD.}}



\maketitle

\begin{abstract}
This paper derives a symbolic multi-body rigid nonlinear model for a space 
telescope using Stoneking's implementation of 
Kane's method. This symbolic nonlinear model is linearized using Matlab 
symbolic functions {\tt diff} and {\tt inv}
because the analytic linearization is intractable for manual derivation. The linearized
system model is then used to design the controllers using both linear quadratic regulator
(LQR) and robust pole assignment methods. The closed-loop systems for the two designs 
are simulated using both the rigid model as well as a second model containing 
flexible modes. The performances of the two 
designs are compared based on the simulation testing results. Our conclusion is that the
robust pole assignment design offers better performance than that of the LQR 
system in terms of actuator usage and pointing accuracy. However, the LQR 
approach remains an effective first design step that can inform the selection 
of real eigenvalues for robust pole assignment. The proposed
method may be used for the modeling and controller designs for various multi-body systems.
\end{abstract}

{\bf keywords}: Modeling, space telescope control, multibody dynamics, LQR, robust pole assignment.

\section{Introduction}
\begin{figure}[!t]
\centering
\includegraphics[width=3.5in]{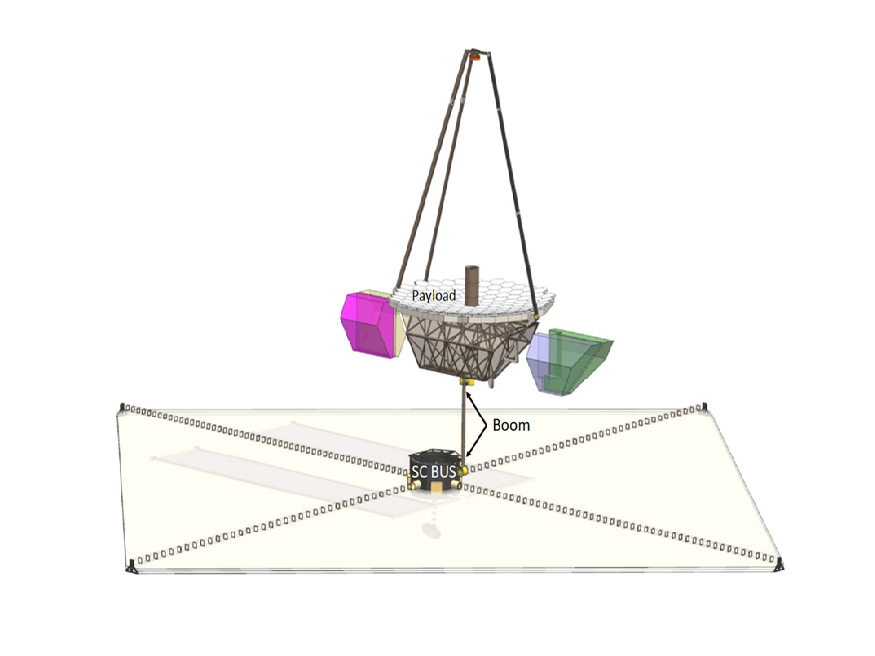}
\caption{The concept of LUVOIR telescope.} \label{conceptD}
\end{figure}

{T}{he} Large UV Optical Infrared Surveyor (LUVOIR) (see Figure \ref{conceptD}), 
to be placed to Sun-Earth L2 point,
is a concept proposed for the key science goal of characterizing a wide
range of exoplanets some of which are potentially habitable. Although the telescope 
is still in the concept phase, NASA has engaged multiple engineering disciplines 
to conduct preliminary design studies \cite{nasa}. The telescope is a typical 
multi-body dynamical system.

Multi-body dynamical systems can be found in many applications including machine design, spacecraft 
dynamics, and robotics. Modern modeling techniques for multi-body dynamics are based on 
d’Alembert's principle in which dynamical systems were essentially converted into static ones through 
the introduction of inertial forces. In 1788, Lagrange formalized this approach by combining the fundamental 
ideas of d’Alembert’s principle with explicit descriptions of virtual work and generalized coordinates \cite{es06}. An extension 
of d’Alembert’s principle valid for holonomic systems was presented in 1909 by Jourdain 
\cite{Jourdain1909}. As many engineering systems are nonholonomic, Kane extended 
d’Alembert’s principle to this general case in 1961 \cite{kane61}. Kane's method has many 
applications particularly in robotics for systems of rigid bodies linked by rotational joints that 
have an arbitrary number of degrees of freedom (see \cite{Buffinton05,cb71,kl85}). Therefore, it is now included in a few 
engineering handbooks, such as \cite{Buffinton05,schiehlen90}. Although Kane's
method has become popular, controversy exists surrounding its originality and efficiency when
compared with the Gibbs-Appell equations \cite{desloge87,lb90}. Recently, Piedboeuf 
indicated that Kane's equations are consistent to the Jourdain principle \cite{piedboeuf93}. 
In \cite{stoneking13}, Stoneking demonstrated that Kane's method can be particularly  
useful in modeling the case of multiple rigid bodies connected via rotary joints, e.g., space telescopes. 
As this was a conference paper, it was limited in both scope as well as exposure. The implementation, 
however, is available as open source software \cite{stoneking10}. 

Using Kane's method as described by Stoneking \cite{stoneking13},
Bentz and Lewis derived a two-body rigid dynamics model for the LUVOIR telescope 
and simulated the initial condition response of a LQR design \cite{bs18}.
This work was followed by similar testing of a higher fidelity three-body rigid dynamics model in \cite{bs18b}. 
However, the derivation is brief and typos exist (which may mislead the readers). In addition,
the full derivation for the three-body model is left to but can be difficult for the readers.
As linked multi-body systems are widely seen in robotics and space applications, we extend the preliminary research 
of \cite{bs18,bs18b} and provide all detailed derivations, which, to our best knowledge, are not
available in any published work. Through the derivation of 
our three-body model, we aim to generate further exposure within the aerospace community 
of Stoneking's implementation of Kane's dynamics and analysis technique that can efficiently model 
rigid multiple bodies, connected via rotary joints having arbitrary degrees of freedom, arranged in tree topologies. 

Although plenty of flexible system modeling methods exist (for example \cite{cp03,lcz20,hg15}),
we are particularly interested in rigid model because the rigid model size is much smaller 
and its states are normally measurable. Therefore, the rigid model is more suitable for the 
control system design than flexible models, and using a rigid model for the controller design 
is widely used in practice. Our ultimate goal is to design a controller for the LUVOIR
telescope in compliance with some arbitrary pointing requirement. As Kane's
multi-body dynamics are nonlinear, and many powerful control techniques such as 
LQR and robust pole assignment are based upon linear models, we have chosen 
to linearize the symbolic model for the purpose of controller design. 
Two controllers are designed based on the linearized model and their 
performances are compared for both rigid and flexible models to give us
some confidence that the designed controller will work for the real system.

There are other multibody modeling methods in the literature. For example,
Li et al \cite{llx22} discussed a flexible multibody spacecraft modeling which 
has a center service module, supporting trusses, and a mirror module. It is assumed
that the center service module and the mirror module are rigid but the trusses 
are flexible. In addition, the rigid center service module's translational motion
is not considered, and connection of the rigid center service module and the
trusses are fixed. Therefore, their model is more specific than the one discussed 
in this paper because we do consider translational motion for all bodies and all
connections are not fixed. Hu et al \cite{hjx12} derived a more general flexible 
multi-body system modeling method, which has much more states. Therefore, 
the model is more suitable for validating the controller design but is not practical 
for controller design.

The rest of the paper is organized as follows: Section \ref{sec:preliminary}
presents some useful background and formulas to be used in the rest of the paper. 
Section \ref{sec:three-body} provides the detailed derivations of the three-body rigid 
nonlinear time invariant dynamics models. Section \ref{sec:design} 
discusses linearization of the dynamics models, and two different controller designs, LQR 
and robust pole assignment, based on the rigid linearized time invariant dynamics model. 
Since the LUVOIR telescope is a flexible structure, we need to validate the designs using 
both the rigid linearized time invariant model as well as the higher fidelity flexible model. 
Redesigns are necessary if a design works for 
simulation testing on the rigid linearized time invariant model but fails the simulation 
testing on the high fidelity flexible model. The final LQR and robust pole 
assignment designs are tested and their performances for the rigid linearized time 
invariant model and flexible model are compared in Section \ref{sec:design} which highlights the advantages of robust pole assignment.

\section{Preliminary}\label{sec:preliminary}

This section provides important concepts and formulas in dynamics theory to be 
used in this paper 	and a brief discussion of Kane's method.

\subsection{Basic concepts and important formulas}

Before we proceed, we present some basic concepts and important formulas which
can be found in \cite{kl85} and will be used repeatedly in the remainder of the paper. 
Let $\B_{\mathcal{F}}=[\b_1,\b_2, \b_3]$ 
be a set of bases of the frame ${\mathcal{F}}$, then a general vector $\vec{\v}$ 
resolved in frame ${\mathcal{F}}$ can be written as:
\begin{equation}
\vec{\v} = v_1\b_1 + v_2 \b_2 +v_3 \b_3
=[\b_1,\b_2, \b_3] [v_1, v_2, v_3]^{\T}=\B_{\mathcal{F}} {\v},
\label{vectorV}
\end{equation}
where $\v=[v_1, v_2, v_3]^{\T}$.
Let $\vec{\omega}_{B/A}$ be the angular rate of frame B relative to frame A resolved
in B. Invoking \cite[(2.3.1)]{kl85}, for any moving vector $\vec{\x}$ resolved in B, its 
derivative in frame A and frame B can be related as:
\begin{equation}
\frac{d \vec{\x}}{d t} \Big|_A = \frac{d \vec{\x}}{d t} \Big|_B 
+ \vec{\omega}_{B/A} \times \vec{\x},
\label{vectorD}
\end{equation}
where $\times$ denotes the cross multiplication of two vectors.
If $\vec{\x}$ is fixed in frame B, then $\frac{d \vec{\x}}{d t} \Big|_B=\0$. Therefore,
we obtain (see \cite[(2.1.2)]{kl85}), 
\begin{equation}
\frac{d \vec{\x}}{d t} \Big|_A = \vec{\omega}_{B/A} \times \vec{\x}.
\label{xFixInB}
\end{equation}
The angular velocity of a rigid body B relative to a reference frame A can be 
expressed in the following form involving $n$ auxiliary references $A_1, \ldots, A_n$
\cite[(2.4.1)]{kl85}:
\begin{equation}
\vec{\omega}_{B/A} = \vec{\omega}_{B/A_1} + \vec{\omega}_{A_1/A_2}
+ \cdots + \vec{\omega}_{A_n/A}.
\label{additionTheorem}
\end{equation}
The angular acceleration of a rigid body B relative to a reference frame A is defined 
as the first time-derivative in A of the angular velocity of $\vec{\omega}_{B/A}$ as 
\cite[(2.5.1)]{kl85}:
\begin{equation}
\vec{\alpha}_{B/A} = \frac{d \vec{\omega}_{B/A}}{dt}.
\label{angularA}
\end{equation}
If $P$ and $Q$ are two points fixed on a rigid body B having an angular velocity
$\vec{\omega}_{B/A}$ relative to a reference frame A, then the velocity of P in A,
denoted as $\vec{\v}_{P/A}$, and the velocity of Q in A, denoted as $\vec{\v}_{Q/A}$,
 are related to each other as \cite[(2.7.1)]{kl85}:
\begin{equation}
\vec{\v}_{P/A} = \vec{\v}_{Q/A} + \vec{\omega}_{B/A} \times \vec{\r},
\label{linearV}
\end{equation}
where $\vec{\r}$ is the position vector from Q to P. The relationship between the acceleration
of P in A, denoted as $\vec{\a}_{P/A}$, and the acceleration of Q in A, denoted as 
$\vec{\a}_{Q/A}$,  is given as \cite[(2.7.2)]{kl85}:
\begin{equation}
\vec{\a}_{P/A} = \vec{\a}_{Q/A} + \vec{\omega}_{B/A}
\times (\vec{\omega}_{B/A} \times \vec{\r}) +\vec{\alpha}_{B/A} \times \vec{\r}.
\label{linearA}
\end{equation}

\subsection{Kane's method}
We will derive the three-body rigid nonlinear model for LUVOIR telescope by using 
Kane's method \cite{stoneking13}. The notations in this section are defined in 
\cite{kl85,stoneking13} and will become clear for the readers who follow the derivation 
to the end of the next section. At that time,  readers will see the beauty of 
Stoneking's form of Kane's method \cite{stoneking13}.
Let $\{ \tau \}$ be the general torque vector of the system, 
$[\J]$ be the general inertia matrix of the system, $\{ \alpha \}$ be the general angular 
acceleration vector of the system, $\{ \omega \}$ be the general angular rate vector of the 
system, $\{ \h \}$ be general angular momentum vector of the system, $\{ \f \}$ be the 
general force vector of the system, $[\M]$ be the general mass matrix of the system, $\{ \a \}$
be the general linear acceleration vector of the system, $\Omega$ be the partial angular velocity
dyad, and $\V$ be the partial velocity dyad ($\Omega$ and $\V$ will be defined
in (\ref{rateMatrix}) and (\ref{velocityMatrix})). 
The Kane's equation in matrix form can be expressed as
\begin{equation}
\Omega^{\T} \left( \{ \tau \} - [\J] \{ \alpha \} - \{ \omega \times \h \} \right)
+\V^{\T} \left( \{ \f \} - [\M] \{ \a \}  \right) = \0,
\label{kaneE}
\end{equation}
where the expression in the first parentheses is Euler's equation, and the expression
in the second parentheses is Newton's second low of motion.
Therefore, formula (\ref{kaneE}) appears at first glance to be trivial; however, there are some significant merits to using Kane's equation for multi-body models as discussed in \cite{stoneking13}. We will see that the 
following relations hold in the rest development.
\begin{subequations}
\begin{gather}
\{ \alpha \} = \Omega \dot{\x}_g + \{ \alpha_r \}, \label{alphaA} \\
\{ \a \} = \V  \dot{\x}_g + \{ \a_r \},   \label{aA}
\end{gather}
\label{accelerations}
\end{subequations}
where $\x_g$ is the generalized speeds of the multi-body system, $\{ \alpha_r \}$ and $\{ \a_r \}$
are items that do not include $\dot{\x}_g$. Substituting equations (\ref{alphaA}) and
(\ref{aA}) into (\ref{kaneE}), and then grouping on $\dot{\x}_g$ yields Stoneking's form of Kane's equation
\begin{eqnarray}
 & & \left(\Omega^{\T} [\J] \Omega + \V^{\T} [ \M ] \V \right) \dot{\x}_g 
\nonumber \\
& = &
\Omega^{\T} \left(   \{ \tau \} - [\J]\{ \alpha_r \} - \{ \omega \times \h \} \right)
\nonumber \\
&  & +\V^{\T} \left( \{ \f \} - [\M] \{ \a_r \}  \right),
\label{kaneS}
\end{eqnarray}
which is the rigid multi-body system model. 
A similar idea was proposed and a similar formula to (\ref{kaneS}) is obtained 
by Hu et al \cite{hjx12} in 2012 for flexible multi-body system modeling. 
The advantages of using Kane's method
for multibody system modeling with tree structure was discussed in
\cite{svl85}. In the next section, we will provide details of using (\ref{kaneS}) for
rigid multi-body system modeling.

\section{Three-body rigid model for LUVOIR telescope}\label{sec:three-body}

\begin{figure}[ht!]
\centering
\includegraphics[width=3.5in]{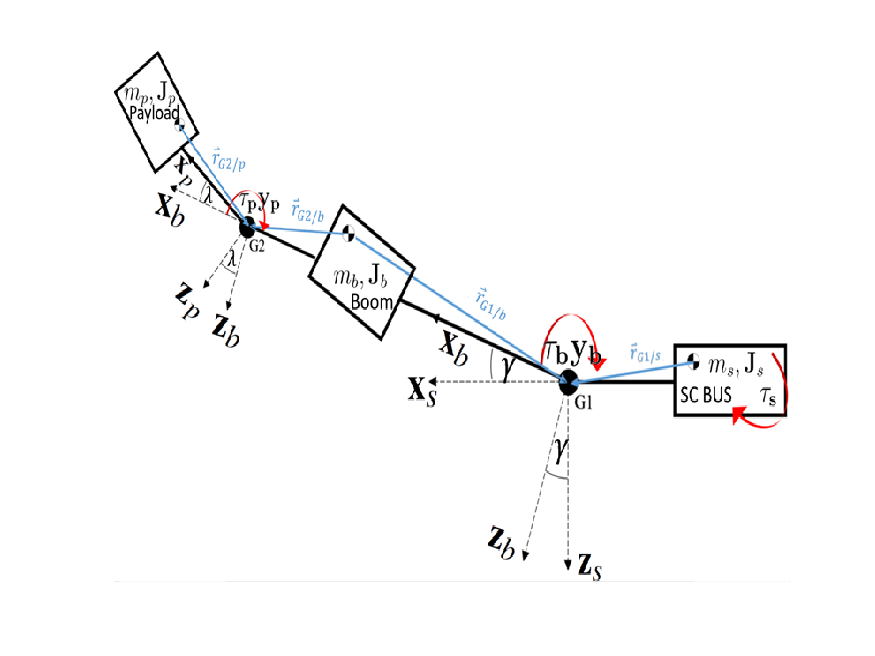}
\caption{The description of the three bodies of LUVOIR telescope.} \label{figureT}
\end{figure}

The LUVOIR-A telescope model is assumed to be composed of three rigid bodies connected 
in serial by two rotary joints as illustrated in Figure \ref{figureT}. 

The three bodies are the spacecraft bus, the boom (tower, or payload articulation system), 
and the payload. Spacecraft bus includes many subsystems such as electrical power system,
propulsion, attitude control system, avionics, command and data handling, thermal 
management system, mechanical and structure. The boom can repoint the payload to any 
position on sky. The payload includes optical telescope assembly, the high definition imager,
the extreme coronagraph for living planetary system, and ultraviolet multi-object spectrograph \cite{nasa}.
Several frames of the LUVOIR telescope will be considered. Let the spacecraft body frame
be denoted as $\mathcal{F}_{s}=[\x_{s}, \y_{s}, \z_{s}]$, the inertial frame be
denoted as $\mathcal{F}_{I}=[\x_{I}, \y_{I}, \z_{I}]$, the boom body frame be
denoted as $\mathcal{F}_{b}=[\x_{b}, \y_{b}, \z_{b}]$, the payload body frame
be denoted as $\mathcal{F}_{p}=[\x_{p}, \y_{p}, \z_{p}]$. The spacecraft frame may
be defined relative to the inertial frame by 
\begin{equation}
\mathcal{F}_{s}^{\T} = \mathcal{O}_{s/I} \mathcal{F}_{I}^{\T}.
\label{S2Iframe}
\end{equation}
where $\mathcal{O}_{s/I}$ is the orientation matrix whose subscript ${s/I}$ represents 
that the orientation of $\mathcal{F}_{s}$ is relative to $\mathcal{F}_{I}$. Using standard 
$3-2-1$ sequence of the intrinsic Euler angle rotations by yaw angle $\psi$, pitch angle 
$\theta$, and roll angle $\phi$, the orientation matrix $\mathcal{O}_{s/I}$ can 
be expressed as
\begin{equation}
\mathcal{O}_{s/I} = \left[ \begin{array}{ccc}
1 & 0 & 0 \\ 0 & \cos(\phi) & \sin(\phi) \\ 0 & -\sin(\phi) & \cos(\phi)
\end{array} \right]
 \left[ \begin{array}{ccc}
\cos(\theta) & 0 & -\sin(\theta) \\ 0 & 1 & 0 \\ \sin(\theta) & 0 & \cos(\theta)
\end{array} \right]
 \left[ \begin{array}{ccc}
\cos(\psi) & \sin(\psi) & 0 \\ -\sin(\psi) & \cos(\psi) & 0 \\ 0 & 0 & 1
\end{array} \right].
\label{S2Irotation}
\end{equation}
Since orientation matrix $\mathcal{O}_{s/I}$ is an orthogonal matrix, we have
\begin{equation}
\mathcal{O}_{I/s} = \mathcal{O}_{s/I}^{\T}
\label{I2Srotation}
\end{equation}

Let the boom gimbal angle be $\gamma$ and the payload gimbal angle be  $\lambda$.
The boom body frame to spacecraft body frame orientation matrix can be expressed as
\begin{equation}
\mathcal{O}_{b/s} = \left[ \begin{array}{ccc}
\cos(\gamma) & 0 & -\sin(\gamma) \\ 0 & 1 & 0 \\ \sin(\gamma) & 0 & \cos(\gamma)
\end{array} \right].
\label{B2Srotation}
\end{equation}
The payload body frame to boom body frame orientation matrix can be expressed as
\begin{equation}
\mathcal{O}_{p/b} = \left[ \begin{array}{ccc}
\cos(\lambda) & 0 & -\sin(\lambda) 
\\ 0 & 1 & 0 \\ \sin(\lambda) & 0 & \cos(\lambda)
\end{array} \right].
\label{P2Brotation}
\end{equation}
The payload body frame to spacecraft body frame orientation matrix can be expressed as
\begin{equation}
\mathcal{O}_{p/s} = \left[ \begin{array}{ccc}
\cos(\gamma+\lambda) & 0 & -\sin(\gamma+\lambda) 
\\ 0 & 1 & 0 \\ \sin(\gamma+\lambda) & 0 & \cos(\gamma+\lambda)
\end{array} \right].
\label{P2Srotation}
\end{equation}

Let the angular velocity of the boom relative to the inertial frame
be denoted as $\vec{\omega}_{b/I}$. We will use similar notations
in the remainder of this paper, for example, $\vec{\omega}_{b/s}$, 
$\vec{\omega}_{p/b}$, and $\vec{\omega}_{s/I}$.
Let $\Gamma_1=[0, 1, 0]^{\T}$ and $\Gamma_2=[0, 1, 0]^{\T}$, $\sigma_1=\dot{\gamma}$ 
and $\sigma_2=\dot{\lambda}$ be the generalized speeds of the rotary joints of
$G1$ and $G2$. Then the angular rate of the rotary joint $G1$ represented 
in the boom frame and the angular rate of the rotary joint $G2$ resolved in the 
payload frame can be written as \footnote{For LUVOIR-B where the connection
between payload and boom has two degrees of freedom, the following equations
may be replaced by (24) in \cite{stoneking13}, but the rest derivation remains essentially the same.}
\begin{equation}
\vec{\omega}_{b/s} =  \vec{\Gamma}_1 \sigma_1, \hspace{0.1in}
\vec{\omega}_{p/b} =  \vec{\Gamma}_2 \sigma_2.
\label{jointRate}
\end{equation}
Using these notations and (\ref{additionTheorem}) (as consistent with \cite[(2.4.1)]{kl85}), 
we have
\begin{equation}
\vec{\omega}_{b/I} =\vec{\omega}_{b/s} +  \vec{\omega}_{s/I}
= \vec{\Gamma}_1 \sigma_1 + \vec{\omega}_{s/I}.
\label{bIRate}
\end{equation}
Let $\B_I$ be the bases of inertial frame, $\B_s$ be the bases of the spacecraft frame,
$\B_b$ be the bases of the boom frame, $\B_p$ be the bases of the payload frame. Then,
we may write (\ref{bIRate}) as
\begin{equation}
\B_b {\omega}_{b/I} 
= \B_b \Gamma_1 \sigma_1  + \B_s {\omega}_{s/I}.
\label{bIRate1}
\end{equation}
By premultiplying $\B_b^{\T}$, we may clear the base dyads and obtain
\begin{equation}
{\omega}_{b/I} 
= \Gamma_1 \sigma_1+ \mathcal{O}_{b/s} {\omega}_{s/I} .
\label{bIRate2}
\end{equation}
Similarly,
\begin{eqnarray}
\vec{\omega}_{p/I} =\vec{\omega}_{p/b} + \vec{\omega}_{b/s} +  \vec{\omega}_{s/I}
=\vec{\Gamma}_2 \sigma_2 + \vec{\Gamma}_1 \sigma_1 +\vec{\omega}_{s/I}.
\label{pIRate}
\end{eqnarray}
We may write (\ref{pIRate}) as 
\begin{eqnarray}
\B_p {\omega}_{p/I} & =& \B_p {\omega}_{p/b}+ \B_b {\omega}_{b/s} + \B_s {\omega}_{s/I} 
\nonumber \\
 & =&  \B_p \Gamma_2 \sigma_2 + \B_b\Gamma_1 \sigma_1 +\B_s {\omega}_{s/I}, 
\label{pIRate1}
\end{eqnarray}
and we may clear the base dyads by pre-multiplying $\B_p^{\T}$ and obtain
\begin{eqnarray}
{\omega}_{p/I} 
= \Gamma_2 \sigma_2 + \mathcal{O}_{p/b}\Gamma_1 \sigma_1 + \mathcal{O}_{p/s} {\omega}_{s/I}.
\label{pIRate2}
\end{eqnarray}
Combining (\ref{bIRate2}) and (\ref{pIRate2}) yields
\begin{equation}
\left[ \begin{array}{c}
\omega_{s/I}  \\ \omega_{b/I} \\ \omega_{p/I}
\end{array} \right]
= \underbrace{ \left[ \begin{array}{cccc}
\I_3 & \0_{31} & \0_{31} & \0_{33} \\
\mathcal{O}_{b/s} & \Gamma_1 & \0_{31} & \0_{33} \\
\mathcal{O}_{p/s} & \mathcal{O}_{p/b} \Gamma_1 & \Gamma_2 &  \0_{33} 
\end{array} \right] }_{\Omega}
\left[ \begin{array}{c}
\omega_{s/I}  \\ \sigma_1 \\ \sigma_2 \\ \v_{s/I}
\end{array} \right],
\label{rateMatrix}
\end{equation}
where $\I_3$ is the three-dimensional identity matrix, $\0_{31}$ is the $3\times 1$
all zero matrix, $\0_{33}$ is the $3\times 3$ all zero matrix, and 
$\v_{s/I} $ is the velocity vector of the center of mass of the spacecraft
relative to the inertial frame.

Now, we consider the linear velocity of the center of the mass for the boom 
and the linear velocity of the center of the mass for the payload. First, we introduce
a notation. For any vector $\a=[a_1, a_2, a_3]^{\T}$, let a skew symmetric matrix related 
to $\a$ be defined as 
\begin{equation}
\a^{\times} = \left[ \begin{array}{ccc}
0 & -a_3 & a_2 \\ a_3 & 0 & -a_1 \\ -a_2 & a_1 & 0
\end{array} \right].
\label{crossM}
\end{equation}
The cross product of two vectors $\a \times \b$ can be written as a multiplication
of the matrix $\a^{\times}$ and the vector $\b$, i.e., $\a^{\times} \b$.
Let $\vec{\v}_{s/I}$ be the velocity of the center of mass of the spacecraft in 
the inertial frame, $\vec{\v}_{b/I}$ be the velocity of the center of mass 
of the boom in the inertial frame, $\vec{\v}_{G1/I}$ be the velocity 
of $G1$ in the inertial frame (see Figure \ref{figureT}), $\vec{\r}_{G1/s}$ be the
position vector from the center of mass of the spacecraft to the joint $G1$,
$\vec{\r}_{G1/b}$ be the position vector from the center of mass of the 
boom to the joint $G1$, $\vec{\v}_{p/I}$ be the velocity of the center of 
mass for the payload in the inertial frame, $\vec{\v}_{G2/I}$ be the velocity 
of $G2$ in the inertial frame, $\vec{\r}_{G2/p}$ be the position vector from 
the center of mass of the payload to the joint $G2$, $\vec{\r}_{G2/b}$ 
be the position vector from the center of mass of the boom to the joint $G2$. 
Note that all these $\vec{\v}$ and $\vec{\r}$ vectors are in  
the inertial frame (see Figure \ref{figureT}). Since $G1$ is a point on both 
the spacecraft and the boom, from (\ref{linearV}) (see \cite[(2.7.1)]{kl85}), 
we have
\begin{subequations}
\begin{gather}
\vec{\v}_{G1/I} = \vec{\v}_{s/I} +\vec{\omega}_{s/I} \times \vec{\r}_{G1/s},   \label{vsI}  \\
\vec{\v}_{G1/I} = \vec{\v}_{b/I}+\vec{\omega}_{b/I} \times \vec{\r}_{G1/b}.  \label{vbI} 
\end{gather}
\label{vsbI}
\end{subequations}
Substituting (\ref{vsI}) into (\ref{vbI}) and invoking (\ref{bIRate}) yield
\begin{eqnarray}
\vec{\v}_{b/I} &=& \vec{\v}_{s/I} +\vec{\omega}_{s/I} \times \vec{\r}_{G1/s}
- \vec{\omega}_{b/I} \times \vec{\r}_{G1/b}
\nonumber \\
&=&  \vec{\v}_{s/I} +\vec{\omega}_{s/I} \times \vec{\r}_{G1/s}
- (\vec{\omega}_{b/s}+\vec{\omega}_{s/I}) \times \vec{\r}_{G1/b}\nonumber \\
\label{G1}
\end{eqnarray}
We may represent each vector in an appropriate basis and write (\ref{G1}) as
\begin{eqnarray}
\B_I {\v}_{b/I} &=& \B_I {\v}_{s/I} +\B_s {\omega}_{s/I} \times \B_I {\r}_{G1/s}
\nonumber \\
& &  
- (\B_b \Gamma_1 \sigma_1+ \B_s {\omega}_{s/I}) \times \B_I {\r}_{G1/b}
\label{G1B}
\end{eqnarray}
Using the notations that ${\r}_{G1/s}|_I = \B_I {\r}_{G1/s}$ (where ${\r}_{G1/s}|_I$
means that the vector ${\r}_{G1/s}$ is expressed in the inertial frame) and 
${\r}_{G1/b}|_I= \B_I {\r}_{G1/b}$, we may clear the base dyads by
pre-multiplying $\B_I^{\T}$ and obtain
\begin{eqnarray}
{\v}_{b/I} &=&  {\v}_{s/I} - {\r}_{G1/s}|_I \times \mathcal{O}_{I/s}{\omega}_{s/I} 
+ {\r}_{G1/b}|_I  \times \mathcal{O}_{I/b}\Gamma_1 \sigma_1
\nonumber  \\ 
& & +{\r}_{G1/b}|_I \times \mathcal{O}_{I/s}{\omega}_{s/I} 
\nonumber \\
&=&   \v_{s/I}|_I +\underbrace{\left[\r_{G1/b}|_I^{\times}    
-  \r_{G1/s}|_I^{\times} \right]\mathcal{O}_{I/s}}_{\V_{21} \in {\bf R}^{3 \times 3}}
\omega_{s/I} \nonumber  \\ 
& & + \underbrace{\r_{G1/b}|_I^ \times \mathcal{O}_{I/b}\Gamma_1}_{\v_{22} 
\in {\bf R}^{3 \times 1}} \sigma_1 .
\label{vb} 
\end{eqnarray}
Applying the same idea to the joint $G2$ and invoking (\ref{linearV}) (see \cite[(2.7.1)]{kl85}), 
we have
\begin{subequations}
\begin{gather}
\vec{\v}_{G2/I} = \vec{\v}_{b/I} +\vec{\omega}_{b/I} \times \vec{\r}_{G2/b},   \label{vbI2}  \\
\vec{\v}_{G2/I} = \vec{\v}_{p/I}+\vec{\omega}_{p/I} \times \vec{\r}_{G2/p}.  \label{vpI2} 
\end{gather}
\label{vpbI}
\end{subequations}
Substituting (\ref{vbI2}) into (\ref{vpI2}) and invoking (\ref{pIRate}) yield
\begin{eqnarray}
\vec{\v}_{p/I} &=& \vec{\v}_{b/I} +\vec{\omega}_{b/I} \times \vec{\r}_{G2/b}
- \vec{\omega}_{p/I} \times \vec{\r}_{G2/p}
\nonumber \\
&=&  \vec{\v}_{b/I} +(\vec{\omega}_{b/s} + \vec{\omega}_{s/I} ) \times \vec{\r}_{G2/b}
\nonumber \\
& &  
- (\vec{\omega}_{p/b}+\vec{\omega}_{b/s}+\vec{\omega}_{s/I}) \times \vec{\r}_{G2/p}
\nonumber \\
&=&   \vec{\v}_{b/I} -\vec{\omega}_{p/b}\times \vec{\r}_{G2/p}
 + \vec{\omega}_{b/s} \times (\vec{\r}_{G2/b}-\vec{\r}_{G2/p})
\nonumber \\
& &   + \vec{\omega}_{s/I}  \times (\vec{\r}_{G2/b}-\vec{\r}_{G2/p}).
\label{G2}
\end{eqnarray}
We may represent each vector in an appropriate basis and write (\ref{G2}) as
\begin{eqnarray}
\B_I {\v}_{p/I} &=& \B_I {\v}_{b/I} - \B_p \Gamma_2 \sigma_2 \times \B_I {\r}_{G2/p}
\nonumber \\
& &   +\B_b \Gamma_1 \sigma_1  \times (\B_I {\r}_{G2/b} - \B_I {\r}_{G2/p})
\nonumber \\
& & +\B_s {\omega}_{s/I}  \times (\B_I {\r}_{G2/b} -\B_I {\r}_{G2/p})
\label{G2B}
\end{eqnarray}
Using the notations that ${\r}_{G2/b}|_I = \B_I {\r}_{G2/b}$ and 
${\r}_{G2/p}|_I= \B_I {\r}_{G2/p}$, and invoking (\ref{vb}), we may clear the base 
dyads by pre-multiplying $\B_I^{\T}$ in (\ref{G2B}) and obtain
\begin{eqnarray}
{\v}_{p/I} &=&  {\v}_{b/I} +{\r}_{G2/p}|_I  \times \mathcal{O}_{I/p}\Gamma_2 \sigma_2
+  ({\r}_{G2/p}|_I - {\r}_{G2/b}|_I ) \times \mathcal{O}_{I/b} \Gamma_1 \sigma_1 
\nonumber \\
& & + ({\r}_{G2/p}|_I - {\r}_{G2/b}|_I ) \times \mathcal{O}_{I/s}{\omega}_{s/I} 
\nonumber \\
&=&  {\v}_{s/I} +\V_{21}{\omega}_{s/I} +\v_{22} \sigma_1 
+{\r}_{G2/p}|_I  \times \mathcal{O}_{I/p}\Gamma_2 \sigma_2
\nonumber \\
& & +  ({\r}_{G2/p}|_I - {\r}_{G2/b}|_I ) \times \mathcal{O}_{I/b} \Gamma_1 \sigma_1 
+ ({\r}_{G2/p}|_I - {\r}_{G2/b}|_I ) \times \mathcal{O}_{I/s}{\omega}_{s/I} 
\nonumber \\
& = &  \v_{s/I} +\underbrace{\left[ \r_{G1/b}|_I^{\times}    
-  \r_{G1/s}|_I^{\times} +\r_{G2/p}|_I^{\times} -\r_{G2/b}|_I^{\times} \right]
\mathcal{O}_{I/s}}_{\V_{31} \in {\bf R}^{3 \times 3}}(\omega_{s/I})
\nonumber \\
& & + \underbrace{\left[ \r_{G1/b}|_I^ \times+ \r_{G2/p}|_I^{\times} 
-\r_{G2/b}|_I^{\times} \right] 
\mathcal{O}_{I/b} \Gamma_1 }_{\v_{32} \in {\bf R}^{3 \times 1}}\sigma_1 
+ \underbrace{\r_{G2/p}|_I^{\times} \mathcal{O}_{I/p}
\Gamma_2}_{\v_{33} \in {\bf R}^{3 \times 1}} \sigma_2.
\label{vp} 
\end{eqnarray}
Combining (\ref{vb}) and (\ref{vp}) yields
\begin{equation}
\left[ \begin{array}{c}
\v_{s/I}  \\ \v_{b/I}  \\ \v_{p/I} 
\end{array} \right]
= \underbrace{ \left[ \begin{array}{cccc}
\0_{33} & \0_{31} & \0_{31} & \I_{3} \\
\V_{21} & \v_{22} & \0_{31} & \I_{3} \\
\V_{31} & \v_{32}  & \v_{33} &  \I_{3} 
\end{array} \right] }_{\V}
\left[ \begin{array}{c}
\omega_{s/I}   \\ \sigma_1 \\ \sigma_2 \\ \v_{s/I} 
\end{array} \right].
\label{velocityMatrix}
\end{equation}
In the sequel, we show that $\Omega$, defined in (\ref{rateMatrix}), and $\V$, defined in 
(\ref{velocityMatrix}), are the same ones defined in the differential equations 
(\ref{accelerations}) which will be used to obtain the multi-body system model (\ref{kaneS}). 
Let $\vec{\alpha}_{s/I}$ be the angular acceleration of the center of the mass of 
the spacecraft relative to the inertial frame, $\vec{\alpha}_{b/I}$ be the angular
acceleration of the center of the mass of the boom relative to the inertial frame, 
and $\vec{\alpha}_{p/I}$ be the angular acceleration of 
the center of the mass of the payload relative to the inertial frame, respectively. 
Taking the derivative for (\ref{bIRate}) and invoking (\ref{vectorD}), we have
\begin{eqnarray}
\vec{\alpha}_{b/I} & = & \frac{d \vec{\omega}_{b/I}}{d t} 
=\frac{d (\vec{\omega}_{b/s}+\vec{\omega}_{s/I})}{d t} 
\nonumber \\
& = & \frac{d  \vec{\omega}_{b/s}}{d t}+\frac{d \vec{\omega}_{s/I}}{d t}
\nonumber \\
& = & \vec{\Gamma}_1 \dot{\sigma_1} +\vec{\omega}_{b/I} \times \vec{\Gamma}_1\sigma_1
 + \vec{\alpha}_{s/I}.
\label{angularBv}
\end{eqnarray}
Representing each vector in an appropriate base yields
\begin{eqnarray}
\B_b {\alpha}_{b/I} = \B_b\Gamma_1 \dot{\sigma_1}
+\B_b {\omega}_{b/I} \times \B_b \Gamma_1\sigma_1 + \B_s {\alpha}_{s/I}.
\label{angularB1}
\end{eqnarray}
Pre-multiplying $\B_b^{\T}$ on both sides of (\ref{angularB1}) clears the base dyads 
and yields
\begin{eqnarray}
{\alpha}_{b/I} & = & \Gamma_1 \dot{\sigma_1}+ 
{\omega}_{b/I} \times \Gamma_1\sigma_1 + \mathcal{O}_{b/s} {\alpha}_{s/I}
\nonumber \\
 & = & \Gamma_1 \dot{\sigma_1}+ \mathcal{O}_{b/s} {\alpha}_{s/I}
+ {\alpha}_{b/I}^r,
\label{angularB}
\end{eqnarray}
where 
\begin{equation}
{\alpha}_{b/I}^r={\omega}_{b/I} \times \Gamma_1\sigma_1,
\hspace{0.1in} {\omega}_{b/I} = \mathcal{O}_{b/s} \omega_{s/I}.
\label{alphaBr}
\end{equation}
Taking the derivative for (\ref{pIRate}) and invoking (\ref{vectorD}) yields
\begin{eqnarray}
\vec{\alpha}_{p/I}  & = & \frac{d \vec{\omega}_{p/I}}{d t} 
=\frac{d (\vec{\omega}_{p/b}+\vec{\omega}_{b/s}+\vec{\omega}_{s/I})}{d t} 
\nonumber \\
 & = & \frac{d  \vec{\omega}_{p/b}}{d t}+\frac{d \vec{\omega}_{b/s}}{d t}
+\frac{d \vec{\omega}_{s/I}}{d t}
\nonumber \\
& = & \vec{\Gamma}_2 \dot{\sigma}_2 +\vec{\omega}_{p/I} \times \vec{\Gamma}_2\sigma_2
+ \vec{\Gamma}_1 \dot{\sigma}_1 \nonumber \\
 &  & +\vec{\omega}_{b/I} \times \vec{\Gamma}_1\sigma_1
+\vec{\alpha}_{s/I}.\nonumber \\
\label{angularPv}
\end{eqnarray}
Representing each vector in an appropriate base yields
\begin{eqnarray}
\B_p {\alpha}_{p/I}  & = &  \B_p \Gamma_2 \dot{\sigma}_2 
+ \B_p {\omega}_{p/I} \times \Gamma_2\sigma_2
+ \B_b \Gamma_1 \dot{\sigma}_1 
\nonumber \\
 &  & + \B_b {\omega}_{b/I} \times \Gamma_1\sigma_1
+ \B_s {\alpha}_{s/I}.
\label{angularP1}
\end{eqnarray}
Pre-multiplying $\B_p^{\T}$ on both sides of (\ref{angularP1}) clears the base dyads 
and yields
\begin{eqnarray}
{\alpha}_{p/I} & = &  \Gamma_2 \dot{\sigma}_2 
+{\omega}_{p/I} \times \Gamma_2\sigma_2
+\mathcal{O}_{p/b} \Gamma_1 \dot{\sigma}_1
\nonumber \\
 &  &+\mathcal{O}_{p/b} {\omega}_{b/I} \times \Gamma_1\sigma_1
+ \mathcal{O}_{p/s} \alpha_{s/I}
\nonumber \\
 & = &   \Gamma_2 \dot{\sigma}_2 
 +\mathcal{O}_{p/b} \Gamma_1 \dot{\sigma}_1
+ \mathcal{O}_{p/s} \alpha_{s/I} + {\alpha}_{p/I}^r.
\label{angularP}
\end{eqnarray}
where 
\begin{eqnarray}
{\alpha}_{p/I}^r & = & {\omega}_{p/I} \times \Gamma_2\sigma_2
+\mathcal{O}_{p/b} {\omega}_{b/I} \times \Gamma_1\sigma_1
\nonumber \\
& = & {\omega}_{p/I} \times \Gamma_2\sigma_2
+\mathcal{O}_{p/b} {\alpha}_{b/I}^r,
\hspace{0.1in} {\omega}_{p/I} = \mathcal{O}_{p/b} {\omega}_{b/I}.\nonumber \\
\label{alphaPr}
\end{eqnarray}
Denote $\alpha_1=\dot{\sigma}_1$, $\alpha_2=\dot{\sigma}_2$, and 
\begin{equation}
\dot{\x}_g =[ \dot{\omega}_{s/I},~\dot{\sigma}_1,~\dot{\sigma}_2,~\dot{\v}_{s/I}]^{\T}.
\label{state}
\end{equation}
Combining (\ref{angularB}) and (\ref{angularP}) yields
\begin{eqnarray}
&&\{ \alpha \} :=  \left[ \begin{array}{c}
\alpha_{s/I}  \\ \alpha_{b/I} \\ \alpha_{p/I}
\end{array} \right]
\nonumber \\
& = &  \underbrace{ \left[ \begin{array}{cccc}
\I_3 & \0_{31} & \0_{31} & \0_{33} \\
\mathcal{O}_{b/s} & \Gamma_1 & \0_{31} & \0_{33} \\
\mathcal{O}_{p/s} & \mathcal{O}_{p/b} \Gamma_1 & \Gamma_2 &  \0_{33} 
\end{array} \right] }_{\Omega}
\underbrace{ \left[ \begin{array}{c}
\dot{\omega}_{s/I}  \\ \dot{\sigma}_1 \\ \dot{\sigma}_2 \\ \dot{\v}_{s/I}
\end{array} \right] }_{\dot{\x}_g}
+ \underbrace{ \left[ \begin{array}{c}
\0_{31}  \\ {\alpha}_{b/I}^r \\ {\alpha}_{p/I}^r
\end{array} \right]}_{\{ \alpha_r \} }
\nonumber \\
& = &  \Omega \dot{\x}_g + \{ \alpha_r \},
\label{rateEq}
\end{eqnarray}
which is equivalent to (\ref{alphaA}). We also showed that $\Omega$ defined in (\ref{rateMatrix})
is the same as the one defined in (\ref{rateEq}) or in (\ref{alphaA}). Next, we derive 
equation (\ref{aA}). Let $\vec{\a}_{s/I}$ be the acceleration of the center of the 
mass of the spacecraft relative to the inertial frame, $\vec{\a}_{b/I}$ be the 
acceleration of the center of the mass of the boom relative to the inertial frame, 
$\vec{\a}_{p/I}$ be the acceleration of the center of the mass of the 
payload relative to the inertial frame, $\vec{\a}_{G1/I}$ be the acceleration of 
the joint $G1$ relative to the inertial frame, and $\vec{\a}_{G2/I}$ be the acceleration of 
the joint $G2$ relative to the inertial frame, respectively. Since $G1$ is a point on both 
the spacecraft and the boom, applying (\ref{linearA}) to the joint $G1$, we have
\begin{subequations}
\begin{gather}
\vec{\a}_{G1/I} = \vec{\a}_{s/I} + \vec{\omega}_{s/I}
\times (\vec{\omega}_{s/I} \times \vec{\r}_{G1/s}) 
+\vec{\alpha}_{s/I} \times \vec{\r}_{G1/s},
\label{aG1s} \\
\vec{\a}_{G1/I} = \vec{\a}_{b/I} + \vec{\omega}_{b/I}
\times (\vec{\omega}_{b/I} \times \vec{\r}_{G1/b}) 
+\vec{\alpha}_{b/I} \times \vec{\r}_{G1/b}.
\label{aG1b}
\end{gather}
\end{subequations}
Substituting (\ref{aG1s}) into (\ref{aG1b}) yields
\begin{eqnarray}
\vec{\a}_{b/I} & = & \vec{\a}_{s/I} + \vec{\omega}_{s/I}
\times (\vec{\omega}_{s/I} \times \vec{\r}_{G1/s}) 
+\vec{\alpha}_{s/I} \times \vec{\r}_{G1/s}
\nonumber \\
 &  &-\vec{\omega}_{b/I}\times (\vec{\omega}_{b/I} \times \vec{\r}_{G1/b}) 
- \vec{\alpha}_{b/I} \times \vec{\r}_{G1/b}.
\label{accelerationBv}
\end{eqnarray}
Representing each vector in an appropriate base and invoking (\ref{angularB1}) yields
\begin{eqnarray}
&& \B_I  {\a}_{b/I} \nonumber \\
 & = & \B_I {\a}_{s/I} +  \B_s {\omega}_{s/I}
\times (\B_s  {\omega}_{s/I} \times \B_I  {\r}_{G1/s}) 
\nonumber \\
 &  &+ \B_s {\alpha}_{s/I} \times \B_I {\r}_{G1/s}
\nonumber \\
& & -\B_b {\omega}_{b/I}\times (\B_b {\omega}_{b/I} 
\times \B_I {\r}_{G1/b})\nonumber \\
 &  &-\B_b {\alpha}_{b/I} \times \B_I {\r}_{G1/b}
\nonumber \\
& = & \B_I {\a}_{s/I} +  \B_s {\omega}_{s/I}
\times (\B_s {\omega}_{s/I} \times \B_I {\r}_{G1/s}) 
\nonumber \\
 &  &+ \B_s {\alpha}_{s/I} \times \B_I {\r}_{G1/s}
\nonumber \\
& & -\B_b {\omega}_{b/I}\times (\B_b {\omega}_{b/I} \times \B_I {\r}_{G1/b})
\nonumber \\
& & -(\B_b\Gamma_1 \dot{\sigma_1}+\B_b {\omega}_{b/I} \times 
\Gamma_1\sigma_1 + \B_s {\alpha}_{s/I})\times \B_I {\r}_{G1/b}.\nonumber \\
\label{accelerationB1}
\end{eqnarray}
Pre-multiplying $\B_I^{\T}$ on both sides of (\ref{accelerationB1}) clears the base dyads 
and yields
\begin{eqnarray}
&& {\a}_{b/I} \nonumber \\
 &  = & {\a}_{s/I}
+\mathcal{O}_{I/s} {\omega}_{s/I} \times ({\omega}_{s/I} \times {\r}_{G1/s}|_I)
\nonumber \\
& & +\mathcal{O}_{I/s}{\alpha}_{s/I} \times {\r}_{G1/s}|_I
\nonumber \\
& & -\mathcal{O}_{I/b} {\omega}_{b/I}\times ( {\omega}_{b/I} 
\times {\r}_{G1/b}|_I)
\nonumber \\
& & -(\mathcal{O}_{I/b} \Gamma_1 \dot{\sigma_1}+\mathcal{O}_{I/b} 
{\omega}_{b/I} \times \Gamma_1\sigma_1
+\mathcal{O}_{I/s} {\alpha}_{s/I})\times {\r}_{G1/b}|_I
\nonumber \\
& = &  {\a}_{s/I} +\underbrace{({\r}_{G1/b}|_I-{\r}_{G1/s}|_I) 
\times \mathcal{O}_{I/s}}_{\V_{21} \in {\bf R}^{3 \times 3}} {\alpha}_{s/I}
\nonumber \\
& & +\underbrace{{\r}_{G1/b}|_I \times \mathcal{O}_{I/b} 
\Gamma_1}_{\v_{22} \in {\bf R}^{3 \times1}} \dot{\sigma_1} +{\a}_{b/I}^r,
\label{accelerationB}
\end{eqnarray}
where 
\begin{eqnarray}
 {\a}_{b/I}^r & = & \mathcal{O}_{I/s} {\omega}_{I/s} \times ({\omega}_{s/I} 
\times {\r}_{G1/s}|_I) \nonumber \\
& & - \mathcal{O}_{I/b} {\omega}_{b/I}\times ( {\omega}_{b/I}
 \times {\r}_{G1/b}|_I) 
\nonumber \\
&  &  +{\r}_{G1/b}|_I \times (\mathcal{O}_{I/b} {\omega}_{b/I} 
\times \Gamma_1\sigma_1)
\nonumber \\
& = & \mathcal{O}_{I/s} {\omega}_{I/s}^{\times} ({\omega}_{s/I}^{\times}
 {\r}_{G1/s}|_I) - \mathcal{O}_{I/b} {\omega}_{b/I}^{\times}
 ( {\omega}_{b/I}^{\times} {\r}_{G1/b}|_I) \nonumber \\
& & +{\r}_{G1/b}|_I^{\times} 
\mathcal{O}_{I/b} \alpha_{b/I}^r.
\label{reminderAb}
\end{eqnarray}
Applying (\ref{linearA}) to the joint $G2$, we have
\begin{subequations}
\begin{gather}
\vec{\a}_{G2/I} = \vec{\a}_{b/I} + \vec{\omega}_{b/I}
\times (\vec{\omega}_{b/I} \times \vec{\r}_{G2/b}) 
+\vec{\alpha}_{b/I} \times \vec{\r}_{G2/b},
\label{aG2b} \\
\vec{\a}_{G2/I} = \vec{\a}_{p/I} + \vec{\omega}_{p/I}
\times (\vec{\omega}_{p/I} \times \vec{\r}_{G2/p}) 
+\vec{\alpha}_{p/I} \times \vec{\r}_{G2/p}.
\label{aG2p}
\end{gather}
\end{subequations}
Substituting (\ref{aG2b}) into (\ref{aG2p}) yields
\begin{eqnarray}
\vec{\a}_{p/I} & = & \vec{\a}_{b/I} + \vec{\omega}_{b/I}
\times (\vec{\omega}_{b/I} \times \vec{\r}_{G2/b}) 
+\vec{\alpha}_{b/I} \times \vec{\r}_{G2/b}
\nonumber \\
& & -\vec{\omega}_{p/I}\times (\vec{\omega}_{p/I} \times \vec{\r}_{G2/p}) 
- \vec{\alpha}_{p/I} \times \vec{\r}_{G2/p}.
\label{accelerationPv}
\end{eqnarray}
Representing each vector in an appropriate base yields
\begin{eqnarray}
\B_I {\a}_{p/I} & = & \B_I{\a}_{b/I} +  \B_b{\omega}_{b/I}
\times (\B_b {\omega}_{b/I} \times \B_I {\r}_{G2/b}) 
\nonumber \\
& & + \B_b{\alpha}_{b/I} \times \B_I {\r}_{G2/b}
\nonumber \\
& & -\B_p {\omega}_{p/I}\times (\B_p {\omega}_{p/I} 
\times \B_I {\r}_{G2/p})\nonumber \\
& & -\B_p {\alpha}_{p/I} \times \B_I {\r}_{G2/p}.
\label{accelerationP1}
\end{eqnarray}
Pre-multiplying $\B_I^{\T}$ on both sides of (\ref{accelerationP1}) to clear the base dyads, 
and substituting (\ref{angularB}), (\ref{angularP}), and (\ref{accelerationB}) into 
the formula yields
\begin{eqnarray}
{\a}_{p/I} &  = & {\a}_{b/I}
+\mathcal{O}_{I/b} {\omega}_{b/I} \times ({\omega}_{b/I} \times {\r}_{G2/b}|_I)
+\mathcal{O}_{I/b}{\alpha}_{b/I} \times {\r}_{G2/b}|_I
\nonumber \\
& & -\mathcal{O}_{I/p} {\omega}_{p/I}\times ( {\omega}_{p/I} 
\times {\r}_{G2/p}|_I) -\mathcal{O}_{I/p} {\alpha}_{p/I} \times {\r}_{G2/p}|_I
\nonumber \\
& = &  {\a}_{s/I}+ ({\r}_{G1/b}|_I-{\r}_{G1/s}|_I) \times \mathcal{O}_{I/s} {\alpha}_{s/I}
+{\r}_{G1/b}|_I \times \mathcal{O}_{I/b} \Gamma_1 \dot{\sigma_1} +{\a}_{b/I}^r 
\nonumber \\
& &  +\mathcal{O}_{I/b} {\omega}_{b/I} \times ({\omega}_{b/I} \times {\r}_{G2/b}|_I)
+\mathcal{O}_{I/b}(\Gamma_1 \dot{\sigma_1}+ \mathcal{O}_{b/s} {\alpha}_{s/I}
+ {\alpha}_{b/I}^r) \times {\r}_{G2/b}|_I
\nonumber \\
& & -\mathcal{O}_{I/p} {\omega}_{p/I}\times ( {\omega}_{p/I} 
\times {\r}_{G2/p}|_I)
 -\mathcal{O}_{I/p} (\Gamma_2 \dot{\sigma}_2  
+\mathcal{O}_{p/b} \Gamma_1 \dot{\sigma}_1
+ \mathcal{O}_{p/s} \alpha_{s/I} + {\alpha}_{p/I}^r) \times {\r}_{G2/p}|_I
\nonumber \\
& = &  {\a}_{s/I} +\underbrace{({\r}_{G1/b}|_I -{\r}_{G1/s}|_I
+ {\r}_{G2/p}|_I-{\r}_{G2/b}|_I) \times \mathcal{O}_{I/s}}_{\V_{31} \in {\bf R}^{3\times 3}} 
{\alpha}_{s/I}
\nonumber \\
& & + \underbrace{({\r}_{G1/b}|_I+{\r}_{G2/p}|_I -{\r}_{G2/b}|_I )
\times \mathcal{O}_{I/b} \Gamma_1 }_{\v_{32} \in {\bf R}^{3\times 1}} \dot{\sigma}_1
+ \underbrace{{\r}_{G2/p}|_I  \times  \mathcal{O}_{I/p} \Gamma_2}_{\v_{33} \in {\bf R}^{3\times 1}} 
 \dot{\sigma}_2 + {\a}_{p/I}^r,
\label{accelerationP}
\end{eqnarray}
where 
\begin{eqnarray}
 {\a}_{p/I}^r & = & {\a}_{b/I}^r -{\r}_{G2/b}|_I \times \mathcal{O}_{I/b}{\alpha}_{b/I}^r 
 \nonumber \\
& & +{\r}_{G2/p}|_I  \times  \mathcal{O}_{I/p} {\alpha}_{p/I}^r
\nonumber \\
& & +\mathcal{O}_{I/b} {\omega}_{b/I} \times ({\omega}_{b/I} \times {\r}_{G2/b}|_I)
\nonumber \\
& &  -\mathcal{O}_{I/p} {\omega}_{p/I}\times ( {\omega}_{p/I} \times {\r}_{G2/p}|_I).
\label{reminderAp}
\end{eqnarray}
Combining (\ref{accelerationB}) and (\ref{accelerationP}) yields
\begin{eqnarray}
& & \{ \a \} := \left[ \begin{array}{c}
\a_{s/I}  \\ \a_{b/I}  \\ \a_{p/I} 
\end{array} \right]
\nonumber \\
& = & \underbrace{ \left[ \begin{array}{cccc}
\0_{33} & \0_{31} & \0_{31} & \I_{3} \\
\V_{21} & \v_{22} & \0_{31} & \I_{3} \\
\V_{31} & \v_{32}  & \v_{33} &  \I_{3} 
\end{array} \right] }_{\V}
 \underbrace{ \left[ \begin{array}{c}
\dot{\omega_{s/I}}   \\ \dot{\sigma}_1 \\ \dot{\sigma}_2 \\ \dot{\v}_{s/I} 
\end{array} \right]}_{\dot{\x}_g}
+\underbrace{  \left[ \begin{array}{c}
\0_{31}  \\ {\a}_{b/I}^r \\ {\a}_{p/I}^r
\end{array} \right]}_{\{ \a_r \} } \nonumber \\
& := & \V \dot{\x}_g + \{ \a_r \},
\label{accelerationMatrix}
\end{eqnarray}
which is (\ref{aA}). We also showed that $\V$ defined in (\ref{accelerationMatrix})
is the same as the one defined in (\ref{velocityMatrix}) or in (\ref{aA}). 

Model (\ref{kaneS}) is a very general multi-body rigid system model. For a three-body 
rigid system like LUVOIR, assume that the $3 \times 3$ inertia matrices for the spacecraft,
the boom, and the payload are given as $\J_s$, $\J_b$, and $\J_p$, then we have
\begin{equation}
[\J] = \left[ \begin{array}{ccc}
\J_s & \0_{33} &  \0_{33} \\  \0_{33} & \J_b &  \0_{33} \\  \0_{33} &  \0_{33} &\J_p
\end{array} \right].
\label{J}
\end{equation}
Assume that the masses of the spacecraft, the boom, and the payload are given as 
$m_s$, $m_b$, and $m_p$, then we have
\begin{equation}
[ \M ] = \left[ \begin{array}{ccc}
m_s\I_3 & \0_{33} &  \0_{33} \\  \0_{33} & m_b\I_3 &  \0_{33} \\  \0_{33} &  \0_{33} & m_p\I_3
\end{array} \right].
\label{mass}
\end{equation}
Assume there are no external forces acting on the rigid bodies, then we have $\{ \f \}=\0$. 
Assume that the control torques $\u$ on the spacecraft, the boom, and the payload are 
$\tau_s$, $\tau_b$, and $\tau_p$,  
i.e., 
\begin{equation}
\u=[\tau_s^{\T},~\tau_b^{\T},~\tau_p^{\T}]^{\T},
\end{equation}
then, we have
\begin{equation}
\{ \tau \} = \left[ \begin{array}{c}
\tau_s-\tau_b 
\\  \tau_b -\tau_p 
\\  \tau_p
\end{array} \right].
\label{torque}
\end{equation}
Finally, we can express $\{ \omega \times \h \}$ in terms of the angular rates $\omega_{s/I}$, 
$\omega_{b/s}$, and $\omega_{p/b}$ of the spacecraft, the boom, and the payload as follows:
\begin{equation}
\{ \omega \times \h \} = \left[ \begin{array}{c}
\omega_{s/I} \times \J_s\omega_{s/I}  \\  
\omega_{b/s} \times \J_b\omega_{b/s}  \\  
\omega_{p/b} \times \J_p\omega_{p/b} 
\end{array} \right].
\label{momemtum}
\end{equation}
Let 
\begin{equation}
\L = \left(\Omega^{\T} [\J] \Omega + \V^{\T} [ \M ] \V \right),
\label{L}
\end{equation}
\begin{equation}
\r_1 = \Omega^{\T} \left(   \{ \tau \} - [\J]\{ \alpha_r \} - \{ \omega \times \h \} \right),
\end{equation}
and
\begin{equation}
\r_2 = \V^{\T} \left( \{ \f \} - [\M] \{ \a_r \}  \right).
\label{R}
\end{equation}
Substituting (\ref{J}), (\ref{mass}), (\ref{torque}), (\ref{momemtum}), (\ref{L}), 
(\ref{R}), and $\{ \f \}=\0$ into (\ref{kaneS}), we have the three-body rigid system model:
\begin{equation}
\L \dot{\x}_g = \r_1+\r_2,
\label{3bodyModel}
\end{equation}
or
\begin{equation}
\dot{\x}_g  = \L^{-1}( \r_1+ \r_2).
\label{threeBodyModel}
\end{equation}

Equation (\ref{threeBodyModel}) looks very simple, but it is a nonlinear system because 
$\L$, $\r_1$ and $\r_2$ have nonlinear components of the state and control
variables. We need a linear system model so that we can apply LQR or robust pole
assignment designs. First, we must rescope the model for the purpose of pure attitude control. We take our generalized speeds $\x_g$, discard the $\v_{s/I}$ component which decouples from the attitude states, and add the spacecraft's Euler angles $\phi$, $\theta$, and $\psi$ in order to define our state vector $\x =[\phi, \theta, \psi, \gamma, \lambda, \omega_{s/I}, \sigma_1, \sigma_2]^{\T}$.
The kinematical differential equations associated with these Euler angles in the reference
inertial frame are given as \cite[Page 429, Space-three 1-2-3]{kll93}:
\begin{eqnarray}
\left[ \begin{array}{c}
\dot{\phi} \\ \dot{\theta} \\ \dot{\psi}
\end{array} \right]
=\left[ \begin{array}{c}
\left< \omega_{s/I}, [1,~ \sin(\phi)\tan(\theta),~ \cos(\phi)\tan(\theta) ]^{\T} \right>
\\
\left< \omega_{s/I}, [0,~ \cos(\phi),~ -\sin(\phi)]^{\T} \right>
\\
\left< \omega_{s/I}, [0,~ \sin(\phi)\sec{\theta},~ \cos(\phi)\sec{\theta}]^{\T} \right>
\end{array} \right], \nonumber 
\label{EulerA}
\end{eqnarray}
where $\left< \a , \b \right>$ denotes the inner product of two vectors of $\a$ and $\b$.
Therefore, the revised state space nonlinear system is given as:
\tiny
\begin{eqnarray}
\dot{\x} := \left[ \begin{array}{c}
\dot{\phi} \\ \dot{\theta} \\ \dot{\psi} \\ \dot{\gamma} \\ \dot{\lambda} 
\\  \dot{\omega}_{s/I} \\ \dot{\sigma}_1 \\ \dot{\sigma}_2 \\ 
\end{array} \right]
= \left[ \begin{array}{c}
\left< \omega_{s/I}, [1,~ \sin(\phi)\tan(\theta),~ \cos(\phi)\tan(\theta) ]^{\T} \right>
\\
\left< \omega_{s/I}, [0,~ \cos(\phi),~ -\sin(\phi)]^{\T} \right>
\\
\left< \omega_{s/I}, [0,~ \sin(\phi)\sec{\theta},~ \cos(\phi)\sec{\theta}]^{\T} \right>
\\ 
\sigma_{1}
\\ 
\sigma_{2}
\\ 
\left[ \I_3,~~~ \0_{35} \right] (\L^{-1} (\r_1+\r_2) )
\\ 
\left[ 0,~ 0,~ 0,~ 1,~ 0,~ 0,~ 0,~ 0 \right] (\L^{-1} (\r_1+\r_2) )
\\ 
\left[ 0,~ 0,~ 0,~ 0,~ 1,~ 0,~ 0,~ 0 \right] (\L^{-1} (\r_1+\r_2) )
\end{array} \right] .
\label{newState}
\end{eqnarray}
\normalsize

\begin{remark}
The procedure of the 3-body modeling can easily be applied to any multibody system with
tree structure, and the modeled system has the structure described in \cite{stoneking13}. It is 
also worthwhile to note that the final state space model discarded some states in $\x_g$
and added some states into $\x$, therefore, the dimensions of $\x_g$ and $\x$ are different.
Finally, (\ref{newState}) involves an analytic inverse matrix $\L^{-1}$ (its computation
will be discussed in the next section), and is slightly different from Stoneking's implementation (\ref{kaneS}).
\end{remark}


\section{Linearization and controller design}\label{sec:design}

To use popular controller design methods, we need to have a linearized rigid dynamics model.

\subsection{Linearization}

Now, we linearize the nonlinear system (\ref{newState}) about a desired new
equilibrium state (when this equilibrium state is attained, $\u=\0$)
so that we will have a symbolic linear system. Assume that this equilibrium state is at 
$\x_d=[\phi_d,~\theta_d,~\psi_d,~\gamma_d,~\lambda_d,~0,~0,~0,~0,~0]^{\T}$ and
the control torques are zeros, i.e., $\u=\0$. Therefore,
\begin{equation}
\dot{\x} = \frac{\partial \f(\x,\u)}{\partial \x} 
\Bigg|_{\begin{array}{c} \x=\x_d \\ \u=\0 \end{array}}
(\x-\x_d)+ \frac{\partial \f(\x,\u)}{\partial \u} 
\Bigg|_{\begin{array}{c} \x=\x_d \\ \u=\0 \end{array}} \u.
\label{linearS}
\end{equation}
where $\phi_d=\pi/2$, and $\theta_d = \psi_d =\gamma_d =\lambda_d=0$.

\begin{remark}
In this case, the target equilibrium state is a $90^o$ rotation of the spacecraft
in roll axis from the current state. Our simulation in the next section will show
that the designed controller works in such a large rotational maneuver. In the
next section, we will discuss a method to obtain the analytic formula for (\ref{linearS}).
\end{remark}

\subsection{Symbolic inverse for linearization}

Clearly, it will be very tedious, if it is not impossible, to find the analytic partial derivatives
for (\ref{linearS}), which involves the calculation of the analytic partial derivatives
of $\L^{-1}$. Bentz and Lewis suggested in \cite{bs18}
using Matlab symbolic function `diff' and the symbolic inverse function `inv' 
for matrix $\L$. For this $8 \times 8$ matrix $\L$, even using Matlab symbolic inverse, the
computation is still too complex to handle. Fortunately, we are only interested in the first
five states in (\ref{threeBodyModel}), (see (\ref{state})), which are the last five states in 
(\ref{newState}), We can use the method proposed in \cite{bs18b}. Let
\[
\L=\left[ \begin{array}{cc}
\L_1 & \L_2 \\  \L_2^{\T} & \L_3
\end{array} \right], \hspace{0.1in}
\p=\r_1+\r_2=\left[ \begin{array}{c}
\p_1 \\ \p_2 
\end{array} \right], \hspace{0.1in}
\dot{\x}_g=\left[ \begin{array}{c}
\dot{\x}_{g,1} \\ \dot{\x}_{g,2}
\end{array} \right], 
\]
then, we have
\begin{equation}
\left[ \begin{array}{cc}
\L_1 & \L_2 \\  \L_2^{\T} & \L_3
\end{array} \right]   
\left[ \begin{array}{c}
\dot{\x}_{g,1} \\ \dot{\x}_{g,2}
\end{array} \right]
=\left[ \begin{array}{c}
\p_1 \\ \p_2 
\end{array} \right].
\label{blockE}
\end{equation}
Solving the second equation of (\ref{blockE}) for $\dot{\x}_{g,2}$ gives
\begin{equation}
\dot{\x}_{g,2} = \L_3^{-1}(\p_2-\L_2^{\T}\dot{\x}_{g,1}).
\label{block1}
\end{equation}
Substituting (\ref{block1}) into the first equation of (\ref{blockE}) gives
\begin{equation}
\dot{\x}_{g,1}=\left(\L_1-\L_2\L_3^{-1}\L_2^{\T}\right)^{-1}
\left(\p_1-\L_2\L_3^{-1}\p_2\right),
\label{block2}
\end{equation}
which involves symbolic inverses of a $3 \times 3$ matrix $\L_3^{-1}$ and a
$5 \times 5$ matrix $\left(\L_1-\L_2\L_3^{-1}\L_2^{\T}\right)^{-1}$.

\subsection{Representation of vectors in inertial frame}

All constants in the multi-body model
are provided by mechanical engineers according to the spacecraft designs.
Some constants in the multi-body model are independent to the frames, 
such as mass of spacecraft, mass of payload, etc., but some constants 
are dependent to the frames. Most likely, the distance vectors in a rigid 
body are given in that rigid body frame, but we need to represent these 
distance vectors in the inertial frame in the model (\ref{block2}) 
as discussed in the previous section.

Let $\r_{G1/s}|_s=a_1 \x_s + a_2 \y_s + a_3 \z_s$ be the position vector from
the center of mass of the spacecraft pointing to the joint $G_1$ represented
in the spacecraft frame, 
$\r_{G1/b}|_b=b_1 \x_b + b_2 \y_b + b_3 \z_b$ be the position vector from
the center of mass of the boom pointing to the joint $G_1$ represented
in the boom frame, 
$\r_{G2/b}|_b=c_1 \x_b + c_2 \y_b + c_3 \z_b$ be the position vector from
the center of mass of the boom pointing to the joint $G_2$ represented
in the boom frame, and
$\r_{G2/p}|_p=d_1 \x_p + d_2 \y_p + d_3 \z_p$ be the position vector from
the center of mass of the payload pointing to the joint $G_2$ represented
in the payload frame. Denote $\r_{G1/s}|_I$ be the the position vector from
the center of mass of the spacecraft pointing to the joint $G_1$ represented
in the inertial frame, $\r_{G1/b}|_I$ be the the position vector from
the center of mass of the boom pointing to the joint $G_1$ represented
in the inertial frame,  $\r_{G2/b}|_I$ be the the position vector from
the center of mass of the boom pointing to the joint $G_2$ represented
in the inertial frame, and  $\r_{G2/p}|_I$ be the the position vector from
the center of mass of the payload pointing to the joint $G_2$ represented
in the inertial frame, then using (\ref{S2Irotation}), (\ref{I2Srotation}), 
(\ref{B2Srotation}), (\ref{P2Brotation}),  and (\ref{P2Srotation}), we have
\begin{subequations}
\begin{align}
\r_{G1/s}|_I& = & \mathcal{O}_{I/s}\r_{G1/s}|_s=\mathcal{O}_{I/s}  [a_1, a_2, a_3]^{\T}
\nonumber \\
& = & \mathcal{O}_{s/I}^{\T} [a_1, a_2, a_3]^{\T}, 
\label{positionA} \\
\r_{G1/b}|_I& = & \mathcal{O}_{I/b}  [b_1, b_2, b_3]^{\T}
=\mathcal{O}_{b/I}^{\T} [b_1, b_2, b_3]^{\T}
\nonumber \\
& = & \left( \mathcal{O}_{b/s} \mathcal{O}_{s/I} \right)^{\T} [b_1, b_2, b_3]^{\T}, 
\label{positionB} \\
\r_{G2/b}|_I & = & \mathcal{O}_{I/b}  [c_1, c_2, c_3]^{\T}
=\mathcal{O}_{b/I}^{\T} [c_1, c_2, c_3]^{\T}
\nonumber \\
& = & \left( \mathcal{O}_{b/s} \mathcal{O}_{s/I} \right)^{\T} [c_1, c_2, c_3]^{\T}, 
\label{positionC} \\
\r_{G2/p}|_I & = & \mathcal{O}_{I/p}   [d_1, d_2, d_3]^{\T}
=\mathcal{O}_{p/I}^{\T} [d_1, d_2, d_3]^{\T}
\nonumber \\
& = & \left(\mathcal{O}_{p/s} \mathcal{O}_{s/I} \right)^{\T} [d_1, d_2, d_3]^{\T}, 
\label{positionD} \\
\r_{G1/p}|_I & = & \r_{G2/p}|_I - \r_{G2/b}|_I +\r_{G1/b}|_I,
\label{positionE}  \\
\r_{s/b}|_I & = & \r_{G1/b}|_I - \r_{G1/s}|_I, 
\label{positionF}  \\
\r_{s/p}|_I & = & \r_{G1/p}|_I - \r_{G1/s}|_I.
\label{positionG}  
\end{align}
\label{positionVectors}
\end{subequations}

Using the definition of (\ref{crossM}), we can write
\begin{equation}
\r_{s/b}|_I^{\times} = \left[ \begin{array}{ccc}
0 & -r_{{s/b}_3} & r_{{s/b}_2} \\ r_{{s/b}_3} & 0 & -r_{{s/b}_1} \\ -r_{{s/b}_2} & r_{{s/b}_1} & 0
\end{array} \right],   
\end{equation}
\begin{equation}
\r_{s/p}|_I^{\times} = \left[ \begin{array}{ccc}
0 & -r_{{s/p}_3} & r_{{s/p}_2} \\ r_{{s/p}_3} & 0 & -r_{{s/p}_1} \\ -r_{{s/p}_2} & r_{{s/p}_1} & 0
\end{array} \right],  
\label{skewMatrices}
\end{equation}
and similarly, we can define $\r_{G1/b}|_I^{\times}$, $\r_{G1/p}|_I^{\times}$, and 
$\r_{G2/p}|_I^{\times}$.

Using the parameters of the LUVOIR telescope, we use a Matlab code (which is provided
in \cite{bs18b}) to generate the rigid linearized time invariant model
\begin{equation}
\dot{\x}=\A \x + \B \u
\label{LTImodel}
\end{equation}
with $\A$ and $\B$ given as follows:
\begin{equation}
\A = \left[ \begin{array}{rrrrrrrrrr}
     0  &   0   &  0   &  0  &   0   &  1  &   0  &   0   &  0  &   0\\
     0  &   0   &  0   &  0  &   0   &  0  &   0  &  -1   &  0   &  0\\
     0   &  0  &   0  &   0  &   0  &   0  &   1   &  0  &   0  &   0\\
     0   &  0  &   0  &   0  &   0  &   0  &   0   &  0  &   1  &   0\\
     0   &  0   &  0  &   0  &   0  &   0   &  0  &   0  &   0  &   1\\
     0  &   0 &    0  &   0  &   0  &   0  &   0  &   0 &    0  &   0\\
     0 &    0  &   0   &  0  &   0  &   0  &   0  &   0   &  0  &   0\\
     0  &   0 &    0  &   0 &    0  &   0  &   0   &  0  &   0   &  0\\
     0  &   0  &   0   &  0 &    0  &   0  &   0  &   0  &   0 &    0\\
     0   &  0   &  0   &  0 &    0  &   0   &  0   &  0   &  0  &   0
\end{array} \right], \nonumber
\end{equation}
\tiny
\begin{equation}
\B =10^{-5}  \left[ \begin{array}{rrrrr}
                   0       &            0       &            0     &              0      &             0\\
                   0       &            0       &            0     &              0      &             0\\
                   0       &            0       &            0     &              0      &             0\\
                   0       &            0       &            0     &              0      &             0\\
                   0       &            0       &            0     &              0      &             0\\
   0.027838 &  0.000983 &  0.011040 & -0.001203 & -0.000149\\
   0.000983 &  0.337771 &  0.007301 & -0.473521  & 0.142849\\
   0.011040 &  0.007301 &  0.082504 & -0.009908  & 0.002876\\
  -0.001203&  -0.473521 & -0.009908 &  0.835383 & -0.391328\\
  -0.000149 &  0.142849 &  0.002876 & -0.391328  & 0.331638
\end{array} \right].
\nonumber
\end{equation}
\normalsize

\begin{remark}
The correctness of the rigid linearized model is indirectly validated when the
controller designed by this rigid model stabilizes a separately developed flexible telescope model.
\end{remark}

\subsection{LQR and robust pole assignment designs}\label{sec:design1}

Using the linearized model obtained from Stoneking's form of Kane's method,
we consider two well-known state-space controller design methods: LQR and robust pole 
assignment. For the linearized control system, the LQR design is to find an optimal
feedback gain matrix $\K_{LQR}$ to minimize the following cost function \cite{lvs12}:
\begin{equation}
J=\frac{1}{2} \int_0^{\infty} (\x^{\T}\Q \x+\u^{\T}{\bf R} \u) dt
\end{equation}
under the state space system constraint $\dot{\x}=\A \x + \B \u$ with $\A$ and $\B$ 
being given in the previous section; while the robust pole 
assignment design is to find an optimal feedback gain matrix $\K_{rpa}$ such that (a)
the close-loop eigenvalues of $(\A -\B \K_{rpa})$ are in the desired locations for the rigid 
model, and (b) the sensitivity to the modeling uncertainty (because of using less accurate 
rigid model in controller design) of the close-loop eigenvalues of $(\A -\B \K_{rpa})$ 
is minimized \cite{yang93,ty96}. Let $\Lambda=\diag(\lambda_1, \ldots, \lambda_n)$ and 
$\X=[\x_1, \x_2, \ldots, \x_n]$ be the closed-loop diagonal eigenvalue matrix and the 
corresponding eigenvector matrix of $(\A -\B \K_{rpa})$, the object of the 
robust pole assignment design is to solve the following optimization problem \cite{yang93,ty96}:
\begin{eqnarray}
\min && \frac{1}{2} \det(\X^{H} \X) \nonumber \\
s.t. && (\A -\B \K_{rpa}) \X = \X \Lambda \\
&& \x_i^{H} \x_i =1, \hspace{0.1in} i = 1, \ldots, n, \nonumber
\end{eqnarray}
where the superscript $H$ is used for complex-conjugate transpose, it reduces to a
transpose if all elements of $\Lambda$ are real.
A very efficient algorithm is developed to solve this problem in \cite{ty96}. A Matlab code that 
implements the algorithm of \cite{ty96} is available on the website of \cite{yang23a}.
For more details on robust pole assignment design discussed in this paper, the readers are 
referred to \cite{knv85,ty96}. A concise description of the robust pole assignment is available 
in \cite[Appendix C]{yang19}.

The LQR design has been widely used in aerospace applications because it is considered
a good choice when energy consumption is a major consideration. Pole assignment design is not 
as popular as the LQR design in this case because it is not clear whether the design will consume 
more energy than the LQR design. However, users have noticed that {\bf robust} pole assignment
design \cite{ty96} normally generates a small feedback gain matrix which is a good sign 
of efficient use of energy. There are two other merits associated with the 
robust pole assignment approach. First, the performance of the closed-loop system is robust to 
modeling errors. This is important because the high fidelity model of LUVOIR used in testing will 
include flexible modes that are ignored during control system design due to modeling complexity. Second, robust pole assignment can predict approximations of closed-loop 
system performance characteristics such as settling time, oscillation frequency, etc., by assigning the 
closed-loop poles in desired areas \cite{db08}. For example, to avoid the oscillations
for a second order system, all poles should be assigned to be real according to \cite[Chapter 5]{db08}. 
For higher order systems, the performance is determined by dominate poles which 
are closer to the imaginary axis, therefore, the dominate poles should be assigned to be real. 
LQR cannot do this. Our strategy is to use
LQR approach as an effective first design step that informs the selection of the real 
eigenvalues for robust pole assignment such that these poles are close to the real parts
of the closed-loop eigenvalues of LQR.

\subsection{Simulation testing on rigid model}

\begin{figure}[ht!]
\centering
\includegraphics[width=0.46\textwidth,height=0.2\textheight]{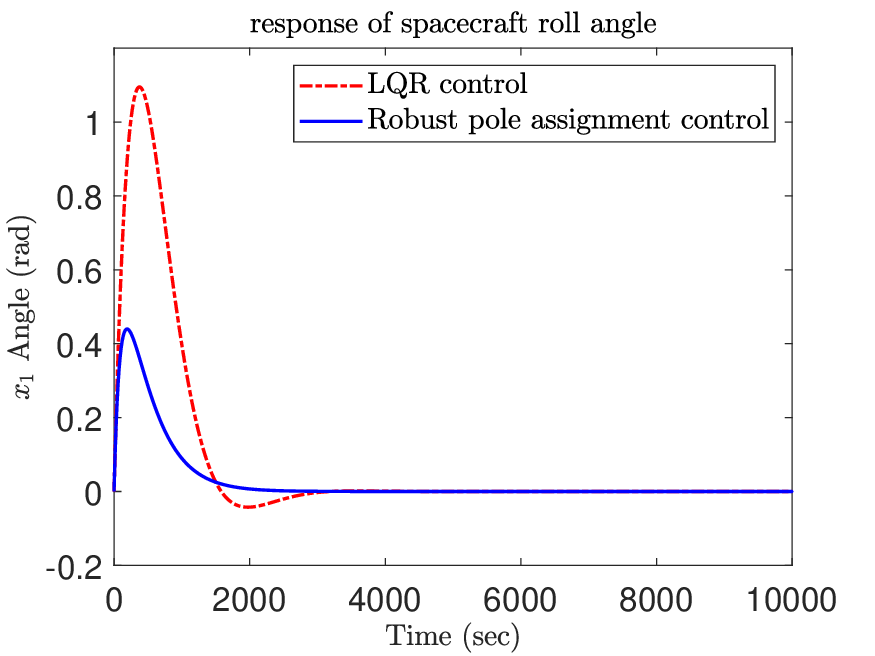}
\hfill
\includegraphics[width=0.46\textwidth,height=0.2\textheight]{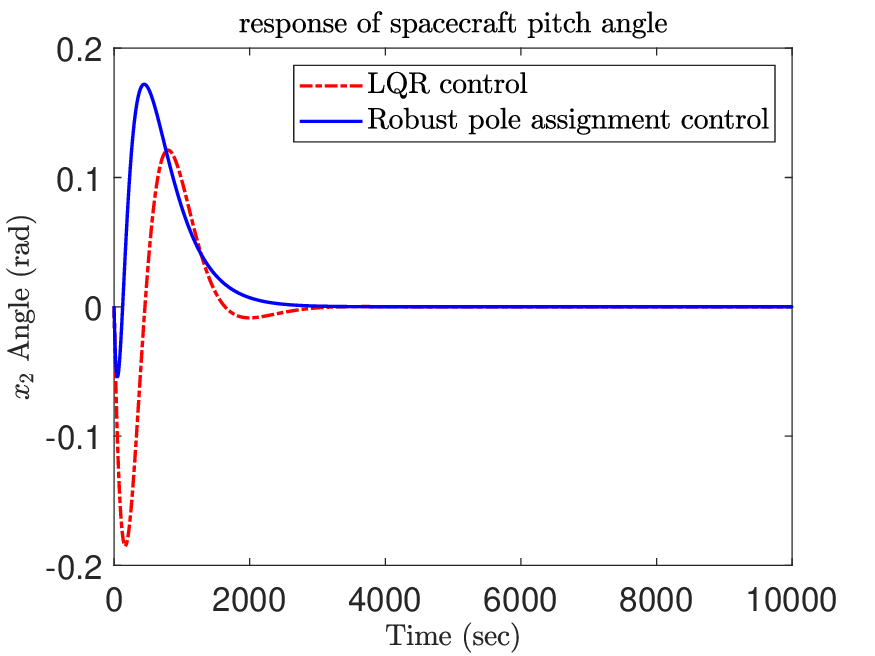}
\caption{LQR and robust pole assignment design comparison for rigid model: 
(a) $x_1$ initial state response (b) $x_2$ initial state response.} \label{figure2}
\end{figure}
\begin{figure}[ht!]
\centering
\includegraphics[width=0.46\textwidth,height=0.2\textheight]{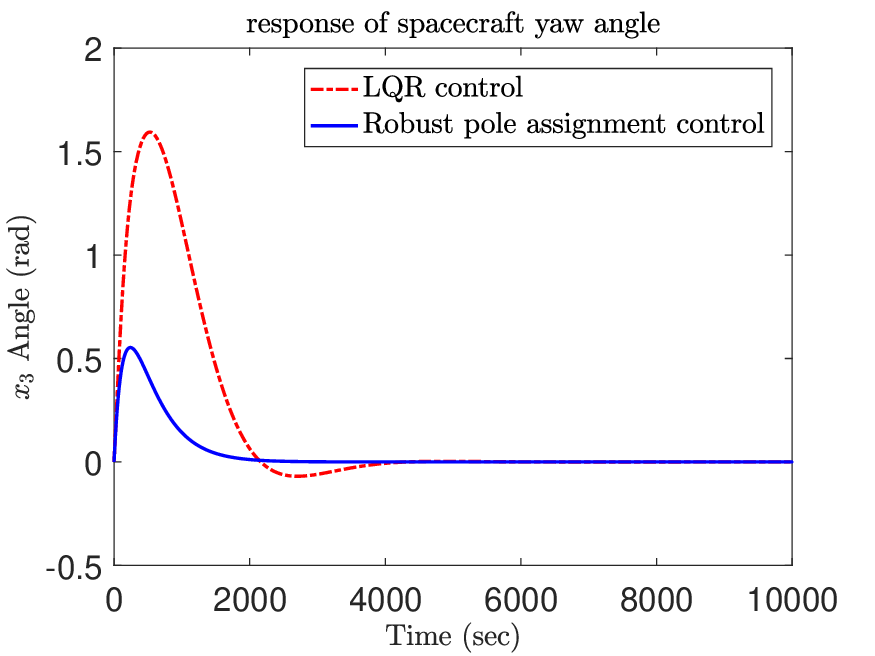}
\hfill
\includegraphics[width=0.46\textwidth,height=0.2\textheight]{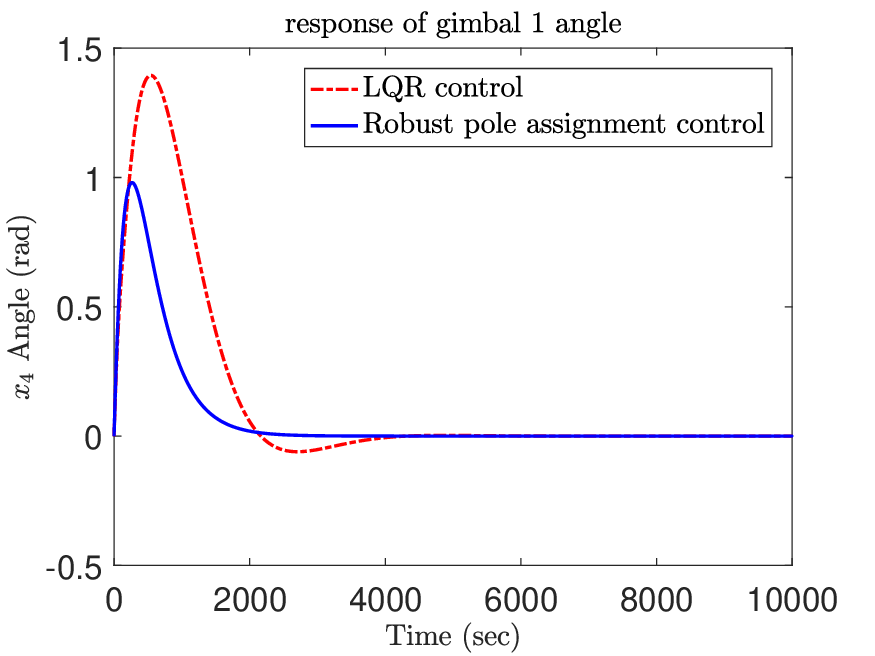}
\caption{LQR and robust pole assignment design comparison for rigid model:
(a) $x_3$ initial state response (b) $x_4$ initial state response.} \label{figure4}
\end{figure}
\begin{figure}[ht!]
\centering
\includegraphics[width=0.46\textwidth,height=0.2\textheight]{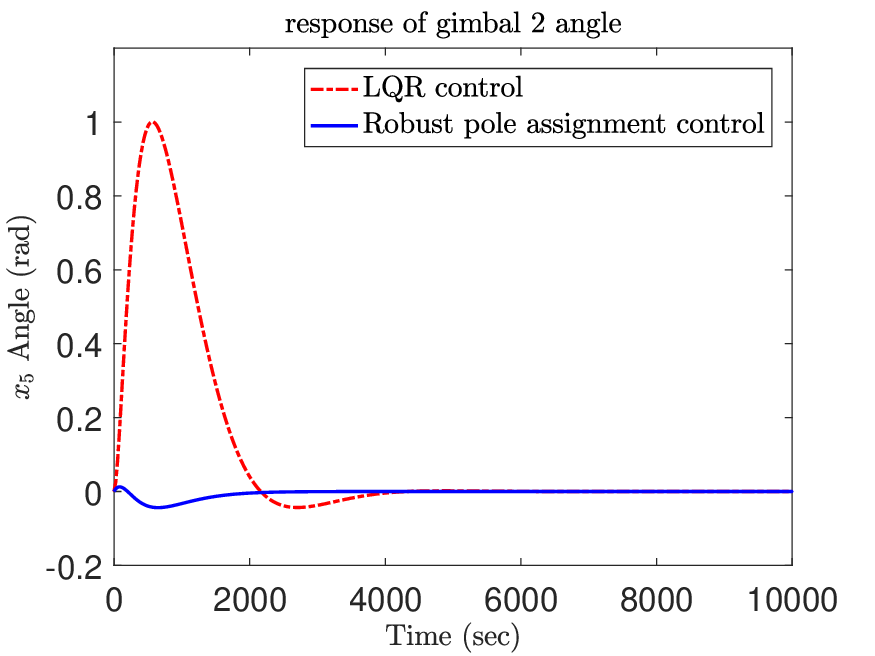}
\hfill
\includegraphics[width=0.46\textwidth,height=0.2\textheight]{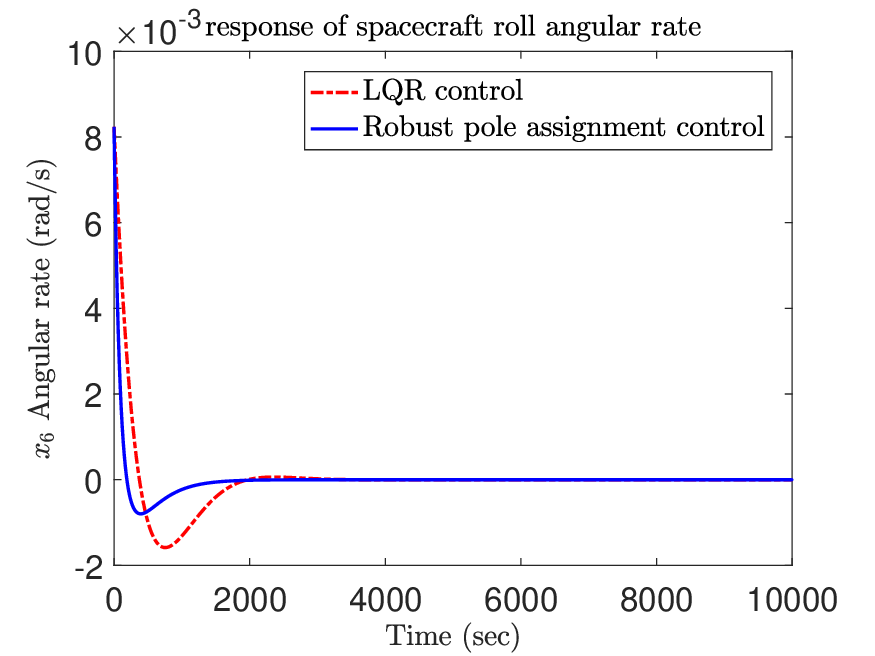}
\caption{LQR and robust pole assignment design comparison for rigid model: 
(a) $x_5$ initial state response (b) $x_6$ initial state response.} \label{figure6}
\end{figure}

We calculated state feedback gain matrices for both LQR and robust pole assignment designs. 
Then, we compared the system performances by analyzing initial state responses and 
energy consumption. The $\Q$ and ${\bf R}$ matrices we used in LQR design are exactly
the same as the ones used in \cite{bs18b}. The target of robust pole assignment is to have a similar settling time to the LQR design but with fewer oscillations\textemdash this is of particular importance to space telescopes with precision pointing requirements. This can be achieved by 
choosing the prescribed closed-loop eigenvalues of the robust pole assignment design 
to have similar real parts to that of the LQR 
design, i.e., placing all closed-loop eigenvalues on the real axis of the complex plane.

For LQR design, the $\Q$ and ${\bf R}$ matrices are selected exactly the same as the
ones in \cite{bs18b}:
\begin{equation}
\Q=\left[ \begin{array} {cc} 1000 \I_8 &  \0_{82} \\ 
\0_{28} & 2000 \I_2 \end{array} \right],
\hspace{0.1in} {\bf R} = \I_5.
\label{QR}
\end{equation}
For robust pole assignment design, the prescribed close-loop eigenvalues are selected as
\footnotesize
\begin{equation}
(-0.0141, -0.0135, -0.0059, -0.0058, -0.0037, -0.0036, -0.0029, -0.0028, -0.00265, -0.00262)
\nonumber
\label{RPA}
\end{equation}
\normalsize
which are close to the real parts of the eigenvalues of $(\A -\B \K_{LQR})$. 
The feedback gain matrices of the LQR and the robust pole assignment are
obtained by a Matlab functions  {\tt lqr} and {\tt robpole} respectively.

\tiny
\begin{eqnarray}
\K_{LQR}=10^4 \left[ \begin{array}{rrrrrrrrrr}
0.0031 &  0.0000  & 0.0000 & -0.0000 &  0.0000 &  1.5433 &  0.0030 & -0.1319 & 0.0045  & 0.0053\\
-0.0000 &  0.0000 &  0.0031 &  -0.0000 &  0.0000 &  0.0030 &  1.2432 & -0.0113 & 0.7782 &  0.4049\\
0.0000 &  -0.0031 &  0.0000 &  -0.0000 &  0.0000 &  -0.1319 &  -0.0113 & 0.8907 &  0.0019 &   0.0000\\
0.0000 &  -0.0000 &   0.0000 &   0.0031 &  -0.0000 &  0.0045 &   0.7782 &   0.0019 &  0.9157 &   0.5169\\
-0.0000 &   0.0000 &  -0.0000 &   0.0000 &  0.0031 &   0.0053&   0.4049 & 0.0000 &  0.5169 &   0.8071 
\end{array} \right] \nonumber
\end{eqnarray}
\begin{eqnarray}
\K_{rpa}=10^4 \left[ \begin{array}{rrrrrrrrrr}
-0.0088 & -0.0013  &  0.0067  & -0.0000 &  -0.0025 &  -3.8298 &   1.8600 &   0.5231 &   0.0202 &  -0.6649\\
0.0052  &  0.0005  & -0.0112  & -0.0025 &  -0.0042 &   1.4826 &  -4.5333 &  -0.0981 &  -1.8008 &  -2.2176\\
0.0013  &  0.0022  & -0.0010  & -0.0002 &   0.0008 &   0.5224 &  -0.2667 &  -1.1273 &  -0.0748 &   0.2379\\
0.0040  &  0.0004  & -0.0089  & -0.0023 &  -0.0044 &   1.1578 &  -3.6498 &  -0.0756 &  -1.6557 &  -2.2215\\
0.0023  &  0.0001  & -0.0056  & -0.0015 &  -0.0045 &   0.6759 &  -2.3400 &  -0.0008 &  -1.1043 &  -2.1248
\end{array} \right] \nonumber
\end{eqnarray}
\normalsize
We demonstrated that the controllers stabilize both the rigid linearized time invariant model
and the flexible LUVOIR telescope model via simulation using Matlab and Simulink.

\subsubsection{Oscillation comparison of the two designs}

\begin{figure}[ht!]
\centering
\includegraphics[width=0.46\textwidth,height=0.2\textheight]{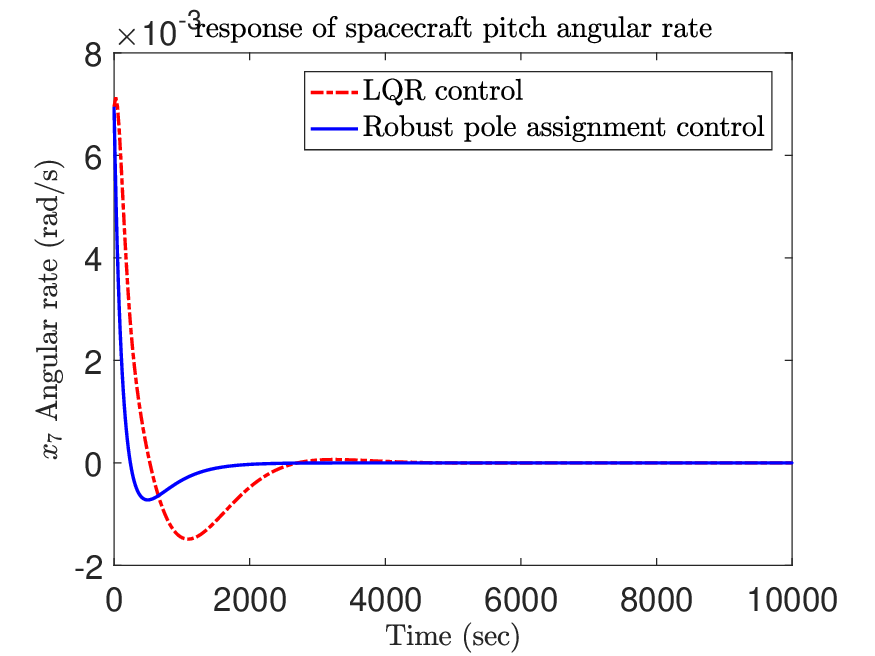}
\hfill
\includegraphics[width=0.46\textwidth,height=0.2\textheight]{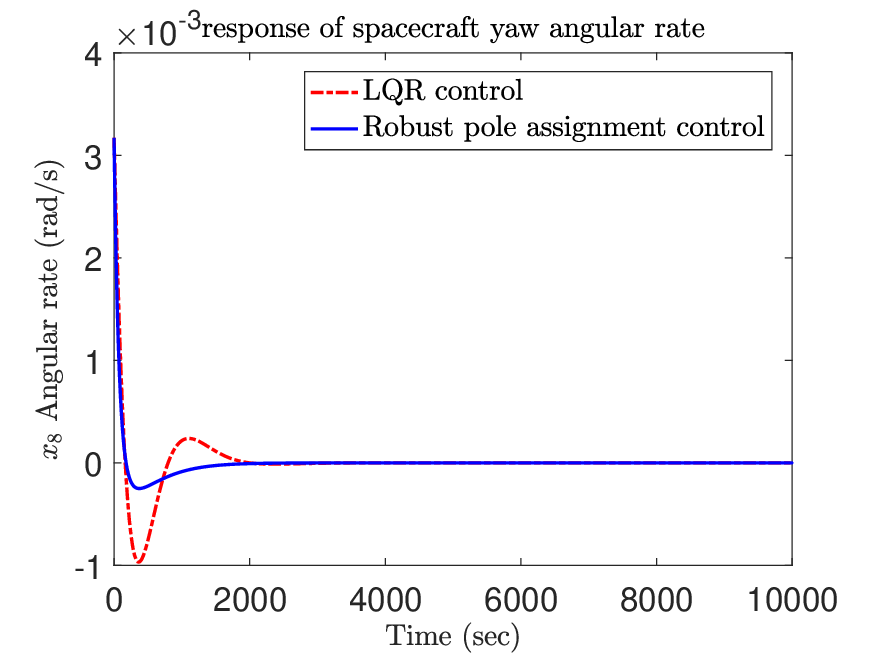}
\caption{LQR and robust pole assignment design comparison for rigid model: 
(a) $x_7$ initial state response (b) $x_8$ initial state response.} \label{figure8}
\end{figure}

\begin{figure}[ht!]
\centering
\includegraphics[width=0.46\textwidth,height=0.2\textheight]{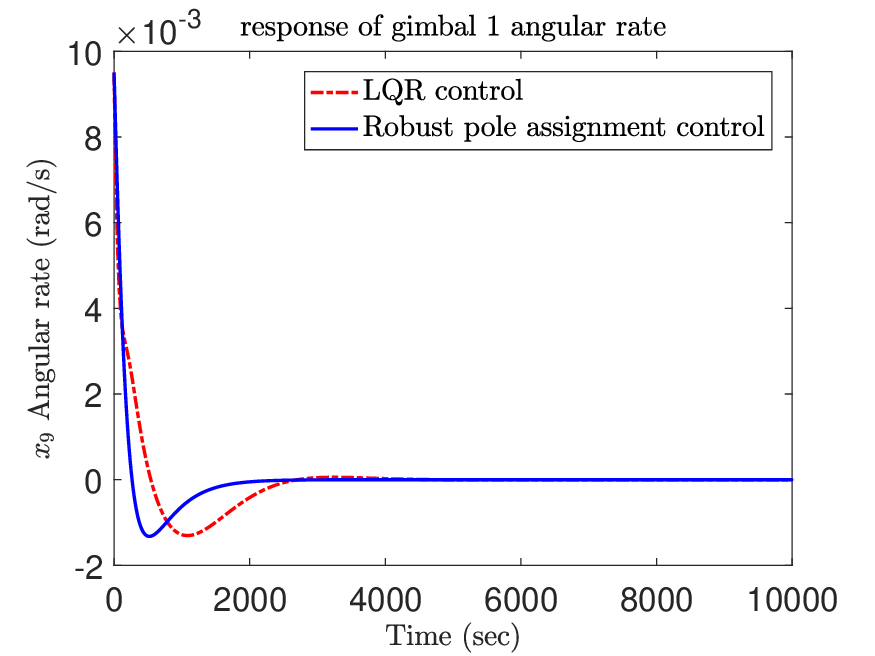}
\hfill
\includegraphics[width=0.46\textwidth,height=0.2\textheight]{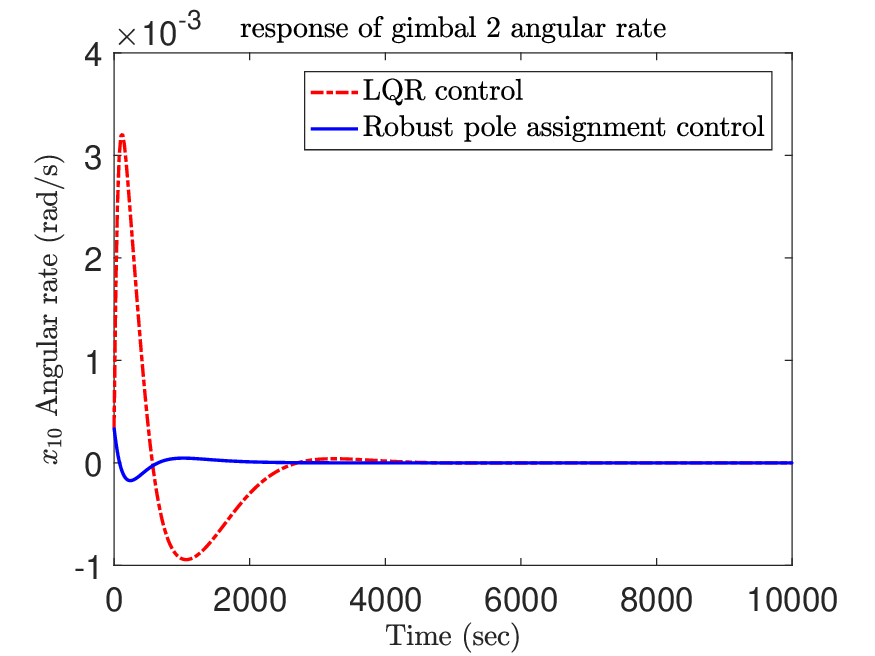}
\caption{LQR and robust pole assignment design comparison for rigid model: 
(a) $x_9$ initial state response (b) $x_{10}$ initial state response.} \label{figure10}
\end{figure}

\begin{figure}[ht!]
\centering
\includegraphics[width=0.46\textwidth,height=0.2\textheight]{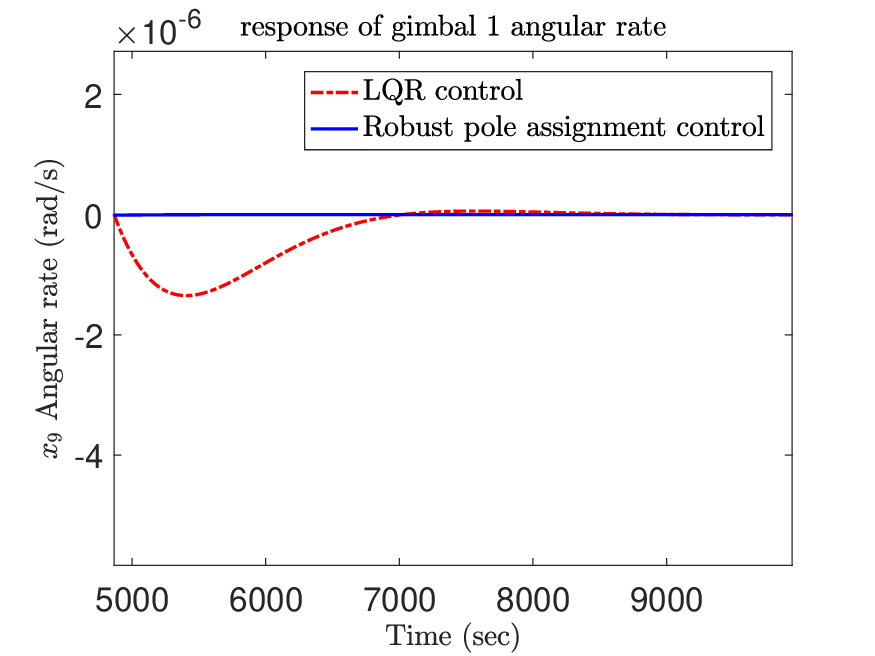}
\hfill
\includegraphics[width=0.46\textwidth,height=0.2\textheight]{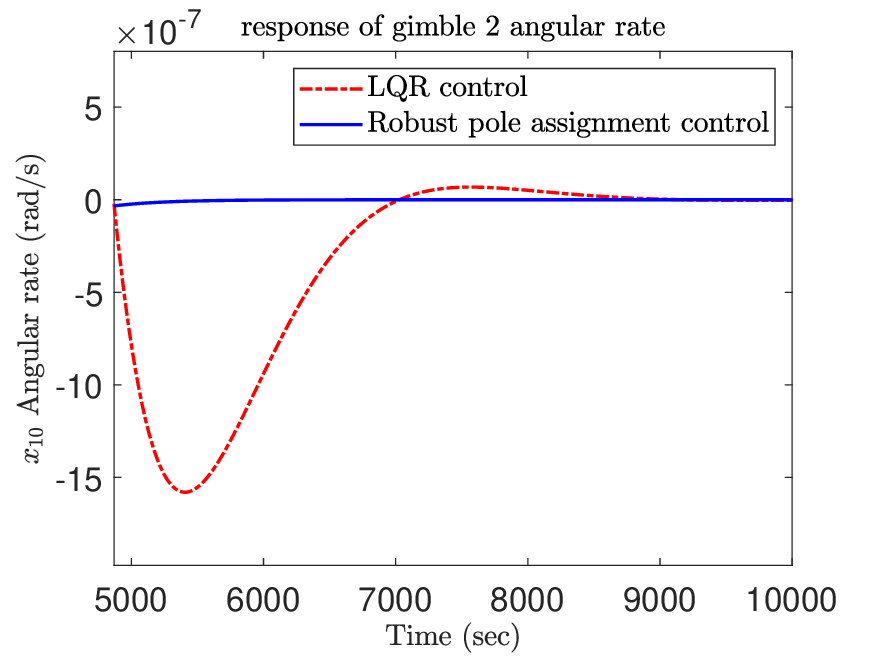}
\caption{LQR and robust pole assignment design comparison for rigid model: (a)  $x_{9}$ initial state response
(b) $x_{10}$ initial state response.} \label{figure12}
\end{figure}

Figures $\ref{figure2}-\ref{figure10}$ compare the three-body rigid linearized model 
initial state responses of the LUVOIR telescope for LQR and robust pole assignment designs. 
It is clear that the initial state responses of robust pole assignment design have fewer 
oscillations in general for all $10$ states (indicating more stable pointing). 
If we amplify the figures, for the LQR design, the oscillation still can be seen after 
$5000$ seconds, but robust pole assignment initial state responses do not have 
this kind of long term oscillations. This meets our expectation as discussed earlier.
Since both controller designs are based on rigid model, when the controllers are applied 
to the rigid model, we don't see the jitter that will be seen when the controllers are
applied to flexible model. Figure \ref{figure12} depicts the long term 
oscillation effects of $x_{9}$ and $x_{10}$. The osculations will adversely affect the 
telescope pointing because disturbances can occur at any time for many different reasons.  

\subsubsection{Energy consumption comparison of the two designs}

Using the least energy consumption to achieve the desired performance is always an 
important design consideration in space missions. Therefore, 
we compare the energy consumption of the two designs. For both LQR and robust pole 
assignment designs, the energy consumption can be measured by
\begin{equation}
\int_0^{\infty} \| \u(t) \|  dt,
\label{fuelConsumption}
\end{equation} 
where $\| \cdot \|$ denotes the Euclidean norm. 
Using formula (\ref{fuelConsumption}), we get 
$\int_0^{10000} \| \u_{LQR} \| dt  =1.45506e+05$ 
and $ \int_0^{10000} \| \u_{rpa} \| dt   =1.28432e+05$, which shows that 
robust pole assignment consumes noticeably less energy.


\subsection{Simulation testing on the flexible model}\label{sec:design2}

\begin{figure}[ht!]
\centering
\includegraphics[width=0.46\textwidth,height=0.2\textheight]{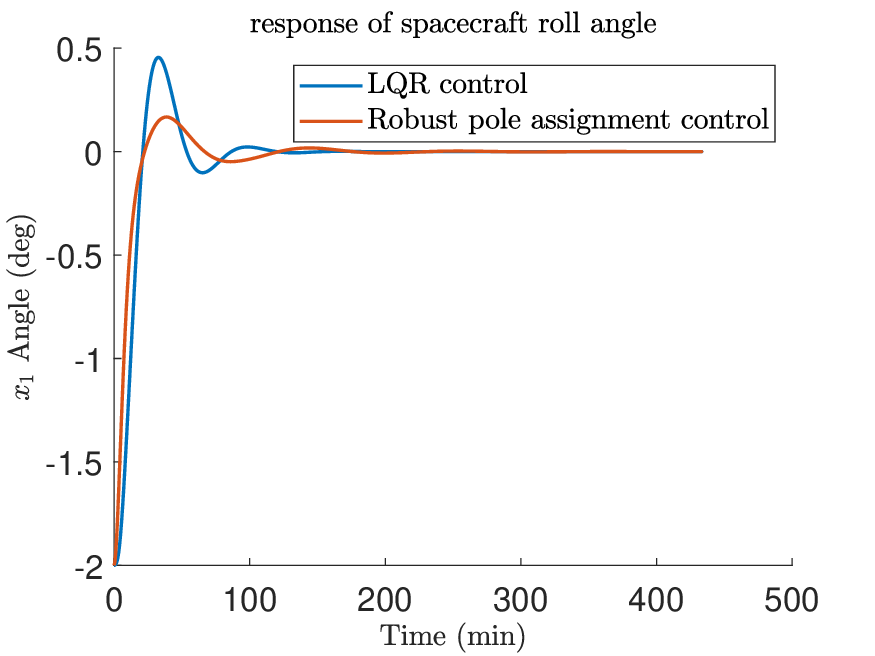}
\hfill
\includegraphics[width=0.46\textwidth,height=0.2\textheight]{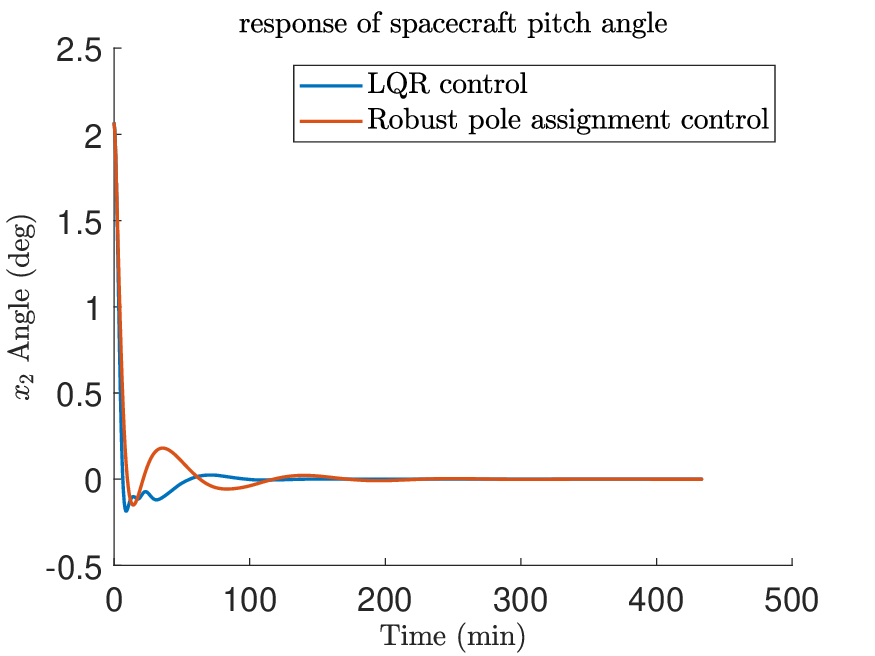}
\caption{LQR and robust pole assignment design comparison for flexible model: 
(a) Roll angle initial state response (b) Pitch angle initial state response.} \label{figure1}
\end{figure}
\begin{figure}[ht!]
\centering
\includegraphics[width=0.46\textwidth,height=0.2\textheight]{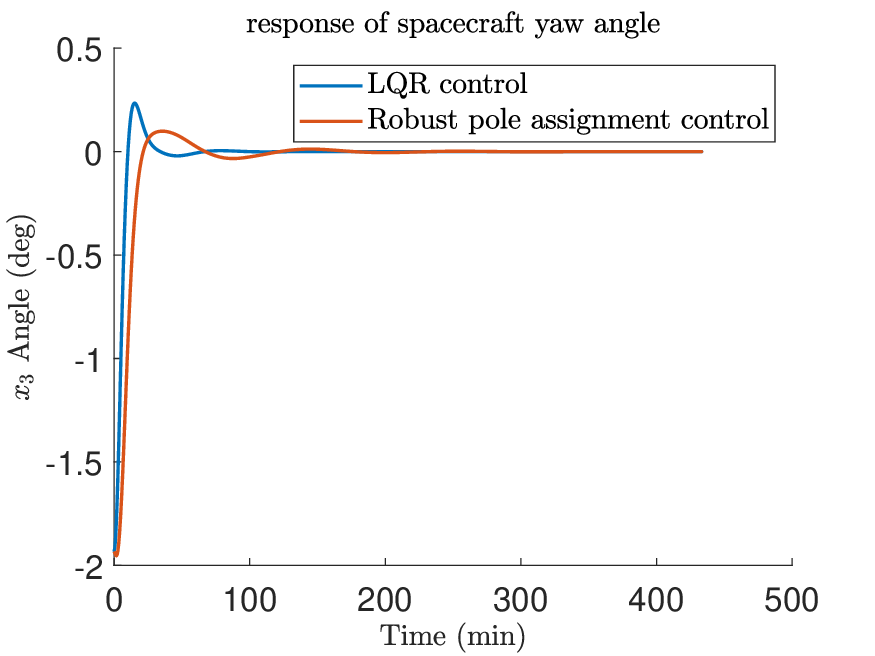}
\hfill
\includegraphics[width=0.46\textwidth,height=0.2\textheight]{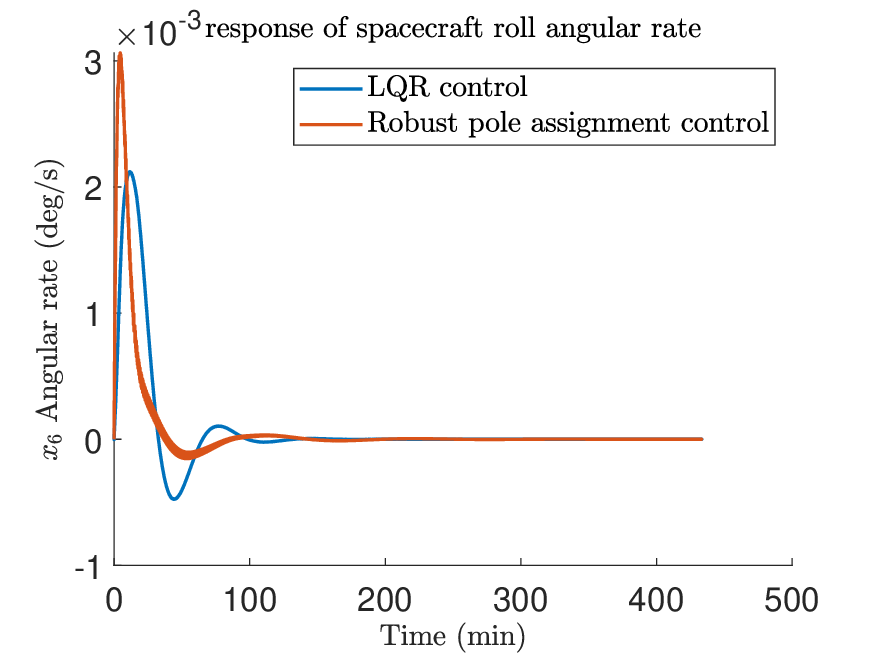}
\caption{LQR and robust pole assignment design comparison for flexible model:
(a) Yaw angle initial state response (b) Roll angular rate initial state response.} \label{figure3}
\end{figure}
\begin{figure}[ht!]
\centering
\includegraphics[width=0.46\textwidth,height=0.2\textheight]{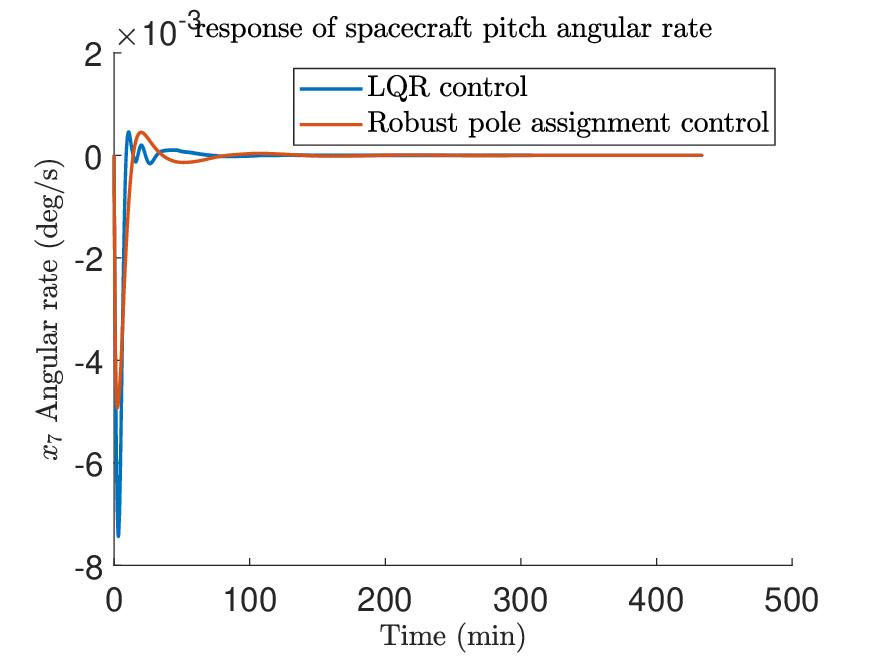}
\hfill
\includegraphics[width=0.46\textwidth,height=0.2\textheight]{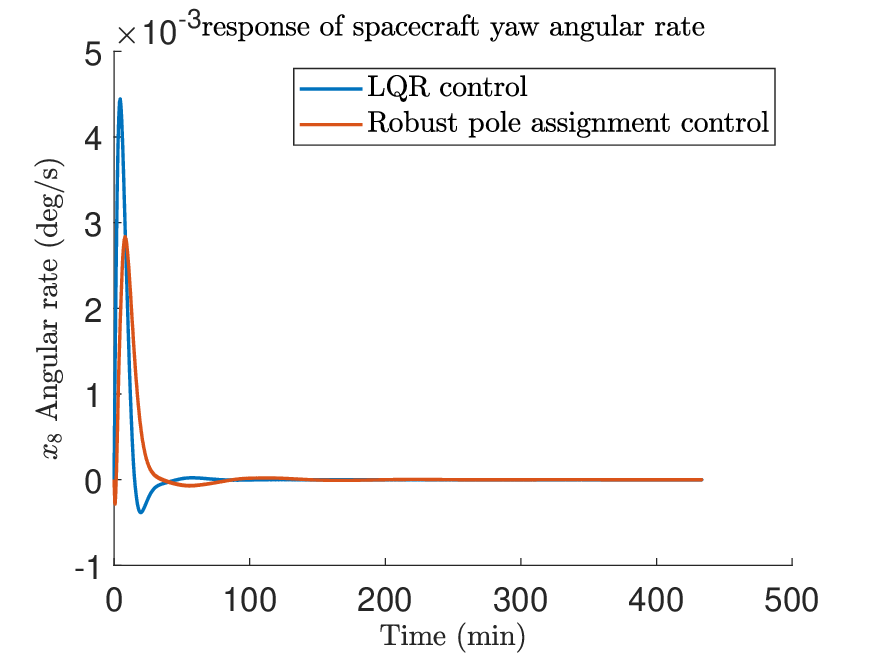}
\caption{LQR and robust pole assignment design comparison for flexible model: 
(a) Pitch angular rate initial state response (b) Yaw angular rate initial state response.} \label{figure5}
\end{figure}
\begin{figure}[ht!]
\centering
\includegraphics[width=0.46\textwidth,height=0.2\textheight]{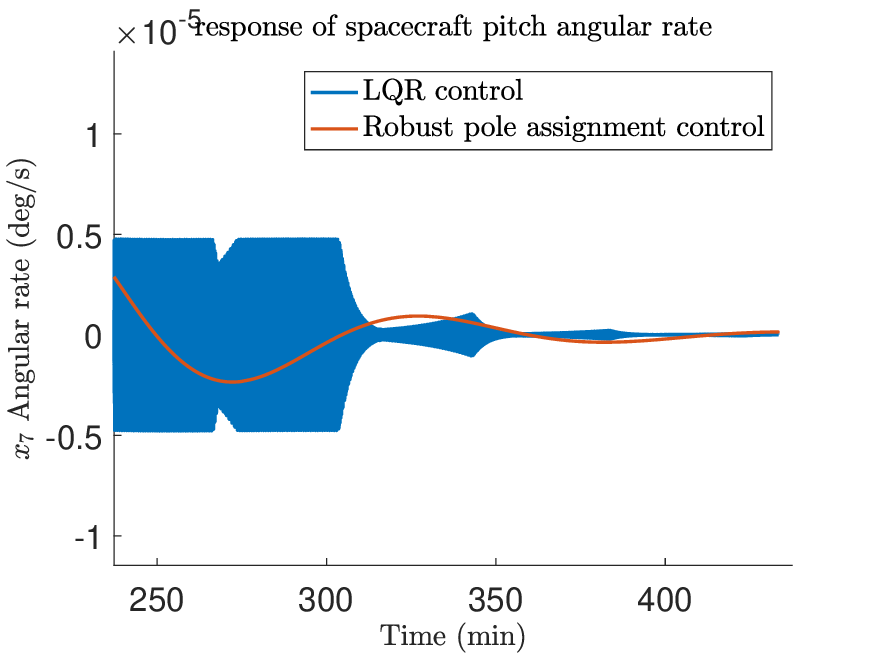}
\hfill
\includegraphics[width=0.46\textwidth,height=0.2\textheight]{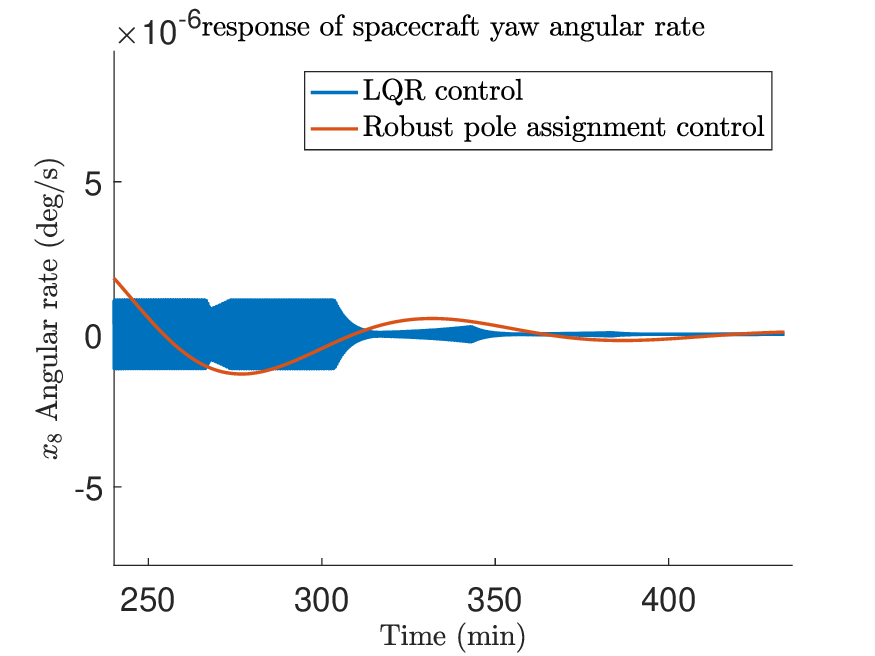}
\caption{LQR and robust pole assignment design comparison for flexible model: 
(a) Pitch augular rate initial state response (b) Yaw angular rate initial state response.} \label{figure7}
\end{figure}
\begin{figure}[ht!]
\centering
\includegraphics[width=0.46\textwidth,height=0.2\textheight]{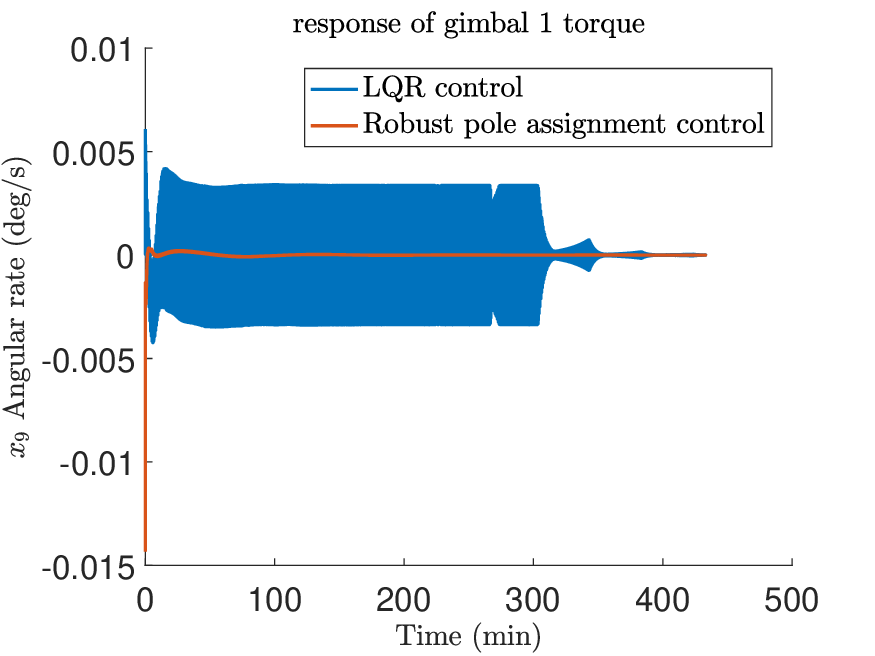}
\hfill
\includegraphics[width=0.46\textwidth,height=0.2\textheight]{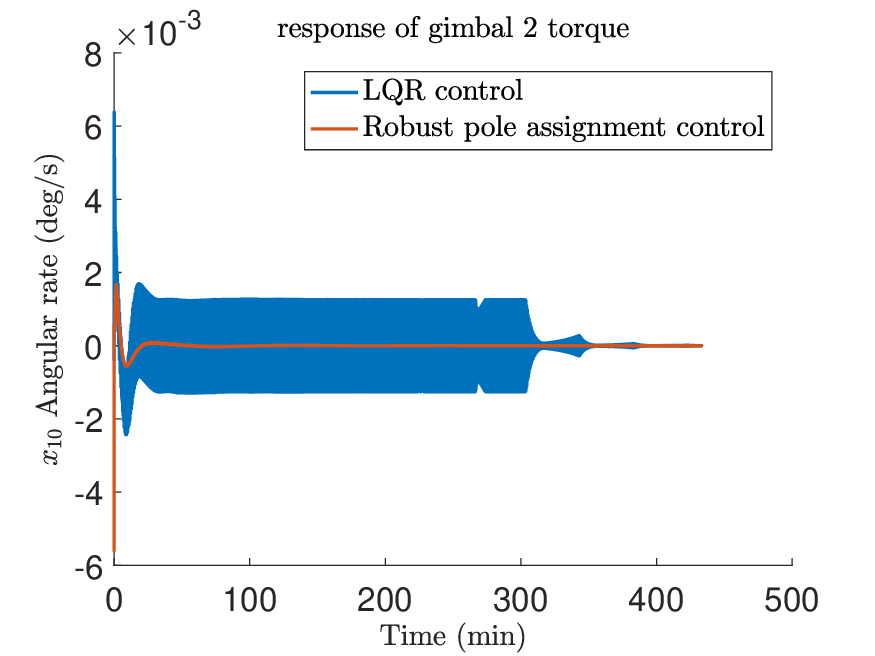}
\caption{LQR and robust pole assignment design comparison for flexible model: 
(a) Gimbal 1 torque initial state response (b) Gimbal 2 torque initial state response.} \label{figure9}
\end{figure}

The LUVOIR flexible model has $20$ modes on the spacecraft, $9$ on the boom, 
and $100$ on the payload \cite{bs18b}, the system has about $260$ states and most of these 
states are not measurable. To have an attainable design and a seamless implementation, 
we designed the controllers based on the coarse rigid model. 
Since the rigid model approximates the flexible model, we need to 
validate the designs by using the simulation for the high fidelity flexible model and examining the 
performances. Modeling a flexible mechanical system has been discussed in \cite{meirovitch01} 
and a Simulink implementation for the LUVOIR telescope was described in \cite{bs18b,chen16}. 
As mentioned in the introduction section, if LQR and/or robust pole assignment designs 
stabilize the rigid model but cannot stabilize the high fidelity flexible model, then, redesigns 
are necessary. As a matter of fact, the feedback gain matrices given in the previous 
section are the ones of the final design which are obtained after a few iterations. 

For the designs given in the previous section, the total energy consumption is 
$1.5448e+03$ for the LQR design and is $1.2683e+03$ for the robust pole assignment 
design. Again, the energy consumption for the robust pole assignment design
is slightly less than the one for the LQR design.

The performances of these two designs are compared and displayed in Figures 
\ref{figure1}-\ref{figure9}. The LQR design settles the spacecraft faster
than the robust pole assignment design, which is good. It can also be seen, from Figures 
\ref{figure9}, that torques requested for 
LQR design have oscillations at about 0.45 Hz, which is not good because
not only does it consume more energy, but it may also introduce jitters to the telescope.
Figures in \ref{figure7} (a) and (b) 
are the amplified pitch and yaw angular rate of
the spacecraft body, which shows the 0.45 HZ oscillations after 200 minutes.

We summarize the systematic space telescope design methodology in the following procedure:
\begin{itemize}
\item[1.] Develop a rigid symbolic nonlinear multibody model using Stoneking’s form of Kane’s equation (\ref{kaneS}).
\item[2.] Take symbolic inverse for Kane's model to obtain the symbolic nonlinear state space model
 $\dot{\x}_g =\f_g(\x, \u)$.
\item[3.] Determine spacecraft kinematical differential equations associated with the 
spacecraft Euler angle using the method provided in \cite{kll93}.
\item[4.] Determine rotary angular dynamics.
\item[5.] Extract relavent states from the symbolic nonlinear state space model $\dot{\x}_g =\f_g(\x, \u)$.
\item[6.] Combine states obtained in Steps 3, 4, and 5 to form a rigid
nonlinear symbolic state space model $\dot{\x} =\f(\x, \u)$. 
\item[7.] Symbolically linearize the nonlinear system about the desired equilibium point to get a 
symbolic rigid linear system model.
\item[8.] Using the spacecraft parameters to populate the symbolic model to get the
 spacecraft specific rigid linear system model.
\item[9.] Design a LQR controller which stabilizes both the rigid linear system model and 
the flexible system model (developed separately from the rigid linear system model).
\item[10.] Design a robust pole assignment controller, whose desired closed-loop poles are all real 
and the value of the real poles are close to the real part of the closed loop poles of LQR design,
such that the robust pole assignment design stabilizes both the rigid linear system model and 
the flexible system model.
\end{itemize}

\section{Conclusion}\label{sec:lastS}
In this paper, we presented a modeling method for multi-body system using
Kane's method. A rigid model for LUVOIR telescope is established as a result.
LQR and robust pole assignment methods are used to design the controllers
using the linearized rigid model. Simulation test of the closed loop system 
using both rigid and flexible models are performed. 
The test result based on the rigid linearized time invariant model shows that 
robust pole assignment has better performance in terms of energy consumption and 
pointing accuracy (measured by the low frequency oscillation around the equilibrium point). 
For the test on the flexible Simulink model (which was developed by Roger 
Chen \cite{chen16}), LQR design has a shorter (better) settling time but the 
gimbal commands have oscillations at about 0.45 Hz which may cause the jitters 
problem and affect the image quality of the telescope; while robust pole
assignment design has a similar gimbal command oscillation problem at the beginning,
it attenuates fast to zero. Overall, we recommend the robust pole assignment technique for this application with the caveat that the LQR approach is effective as a method of jump-starting the robust pole assignment design. That is, the LQR approach allows the designer to make an initial pass in which setting time and other performance characteristics are tuned through the traditional cost function weight matrices. This provides a set of desired real eigenvalue components that can be targeted through robust pole assignment in order to refine the performance, e.g., in order to improve damping. 

\section{Acknowledgments}

The authors are grateful to the anonymous reviewers for their constructive comments which
are invaluable in the preparation of the final version. The first author thanks Eric Stoneking and Gary Mosier of
Goddard Space Flight Center at NASA for introducing him into this interesting project.

%
%
%
%
%

\vfill

\end{document}